\documentclass[reqno, 10pt]{amsart}
\usepackage{amsmath}
\usepackage{amssymb}
\usepackage[all]{xy}
\usepackage{hyperref}
\usepackage[capitalise,noabbrev]{cleveref}
\usepackage{mathtools}
\usepackage{tikz-cd}
\usepackage[shortlabels]{enumitem}

\theoremstyle{plain}
\newtheorem{thm}{Theorem}
\newtheorem{thm*}{Theorem}
\newtheorem{cor}[thm]{Corollary}
\newtheorem{prop}[thm]{Proposition}
\newtheorem{lemma}[thm]{Lemma}

\theoremstyle{definition}
\newtheorem{defn}[thm]{Definition}
\newtheorem{rmk}[thm]{Remark}

\numberwithin{thm}{section}

\usepackage{comment}

\def\del{{\partial}}
\def\maps{{\colon}}

\def\<{\langle}
\def\>{\rangle}

\def\Z{{\mathbb Z}}

\def\cal{\mathcal}

\def\cH{{\cal H}}

\def\cN{{\cal N}}

\def\cS{{\cal S}}

\def\cV{{\cal V}}

\def\im{\operatorname{im}}

\def\Aut{\operatorname{Aut}}
\def\hAut{\operatorname{hAut}}
\def\Homeo{\operatorname{Homeo}}
\def\Diffeo{\operatorname{Diffeo}}

\def\lim{\operatorname{lim}}
\def\rk{\operatorname{rk}}

\def\ev{\operatorname{ev}}
\def\ra{\rightarrow}

\def\PD{\operatorname{PD}}

\def\Graph{\operatorname{Graph}}

\def\Spin{\operatorname {Spin}}
\def\BSpin{\operatorname {BSpin}}

\def\BTopSpin{\operatorname {BTopSpin}}
\def\SO{\operatorname {SO}}
\def\BSO{\operatorname {BSO}}
\def\Top{\operatorname {Top}}
\def\Diff{\operatorname {Diff}}

\def\G{\operatorname {G}}
\def\BTop{\operatorname {BTop}}

\def\BSCAT{\operatorname {BSCAT}}
\def\CAT{\operatorname {CAT}}
\def\SCAT{\operatorname {SCAT}}
\def\CATSpin{\operatorname {CATSpin}}
\def\BCAT{\operatorname {BCAT}}
\def\BCATSpin{\operatorname {BCATSpin}}

\def\BTopSpin{\operatorname {BTopSpin}}
\def\BSO{\operatorname {BSO}}
\def\BSpin{\operatorname {BSpin}}

\def\imm{\operatorname{imm}}
\def\reld{\operatorname{rel\del}}

\def\red{\operatorname{red}}
\def\HE{\operatorname{HE}}
\def\proj{\operatorname{proj}}
\def\geom{\operatorname{geom}}
\def\alg{\operatorname{alg}}

\def\id{\operatorname{id}}

\def\NN{{\mathbb N}}

\def\QQ{{\mathbb Q}}

\def\ZZ{{\mathbb Z}}

\def\fs{{\mathfrak s}}

\def\cal{\mathcal}

\def\cH{{\cal H}}

\def\cN{{\cal N}}

\def\cS{{\cal S}}

\def\cV{{\cal V}}

\hypersetup{colorlinks=true,
    linkcolor=blue,
    citecolor=blue,
    pdftitle={Title},
    pdfauthor={Ansel Goh}}

\begin{document}

\title{$h$-cobordisms between $1$-connected manifolds with boundary}
\author[Goh]{Ansel Goh}

\address{University of Texas at Austin, \\ Austin TX}
\email{anselgoh@utexas.edu}

\begin{abstract}
We classify $h$-cobordisms between compact simply-connected $4$-manifolds with boundary in both the smooth and topological categories. This extends a result of Kreck, which proved the corresponding result for closed manifolds. As an application, we recover the injectivity portion of Orson-Powell's computation of the mapping class group for simply-connected topological $4$-manifolds with boundary. 
\end{abstract}

\maketitle

\section{Introduction} 
In \cite{Kreck_hcob}, Kreck showed that the set of diffeomorphism (resp. homeomorphism) classes rel. boundary of $h$-cobordisms between smooth (resp. topological) simply-connected closed $4$-manifolds is in bijection with the set of isometries of their intersection forms, a result that also follows from Lawson's work on invertible $h$-cobordisms \cite{lawson}. Later work of Hambleton-Kreck in \cite{hambletonkreck} classified $h$-cobordisms between closed manifolds with odd order fundamental group up to extension. We prove the analogous result for simply-connected $4$-manifold with boundary. 

Before we state our result, we loosely define some of the relevant terms, which will be covered more precisely in \cref{invts}. First of all, we define an $h$-cobordism $W$ between simply-connected smooth (resp. topological) $4$-manifolds $M_0,M_1$ to be a \textit{$h$-cobordism rel. $h$}, where $h$ is a diffeomorphism (resp. homeomorphism) from $\del M_0$ to $\del M_1$, if we have a decomposition
\[\del W \cong \del M_0 \cup_{\id} \del M_0 \times [0,1] \cup_h \del M_1.\] We denote the set of smooth (resp. topological) $h$-cobordisms rel. $h$ to be $\cH^{\Diff}(M_0,M_1, h)$ (resp. $\cH^{\Top}(M_0,M_1,h)$).

Now, denote the set of orientation preserving homotopy equivalences $M_0\to M_1$ restricting to $h$ on $\del M_0$ up to homotopy rel. $\del M_0$ to be $\HE(M_0,M_1,h)$. In the case of $\HE(M_0,M_0,\id)$, this has a group structure given by composition and we denote it instead as $\hAut(M_0, \del M_0)$. Let $\CAT\in \{\Top, \Diff\}$. We outline the construction of maps 
\[\cH^{\CAT}(M_0,M_1,h) \to \HE(M_0,M_1,h) \to \hAut (M_0, \del M_0).\]

To start, given an $h$-cobordism rel. $h$, \cref{varphi} constructs a map 
\[\cH^{\CAT}(M_0, M_1, h)\to \HE(M_0, M_1, h).\]
We do this by defining a homotopy equivalence $\varphi_W=i_1^{-1}\circ i_0$, where $i_0$ is the inclusion $M_0\hookrightarrow W$ and $i_1^{-1}$ is a particular choice of homotopy inverse of the inclusion of $i_1\maps M_0\hookrightarrow W$ such that $\varphi_W \vert_{\del M_0}=h$. We then show that this definition is independent of choice of $i_1^{-1}$ up to homotopy rel. boundary, so we can define the map via $W\mapsto \varphi_W$. 

Note that, given an $h$-cobordism rel. $h$ from $M_0$ to $M_1$, the $h$-cobordism theorem (\cite[Theorem 7.1A]{FreedmanQuinn}) provides a homeomorphism $F\maps M_1\to M_0$ restricting to $h^{-1}$ on the boundary. Hence, we obtain a map 
\[\cH^{\CAT}(M_0, M_1, h)\to \hAut(M_0,\del M_0), \qquad W \mapsto F\circ \varphi_W.\] When $M_0 = M_1$ and $h=\id$, we give $\cH^{\CAT}(M_0,M_1, h)$ a group structure via the operation of stacking $h$-cobordisms on top of one another and we obtain a canonical choice for $F$, namely the identity. In this scenario, the map becomes a group homomorphism $\cH^{\CAT}(M_0,M_1, h)\to \hAut (M_0, \del M_0)$. 

Next we define our two invariants. First, given a $4$-manifold $M$, let $f\maps M \to M$ be a continuous map restricting to the identity on the boundary. Then, the \textit{variation} associated to $f$ is the map $\Delta_f\maps H_2(M, \del M)\to H_2(M)$ defined by $[x]\mapsto [x-f(x)]$. When $f$ is a homotopy equivalence, we show that $\Delta_f$ belongs to a group of functions $H_2(M, \del M)\to H_2(M)$ defined by Saeki in \cite{Saeki} called \textit{Poincar\'e variations}, which we denote by $\cV(H_2(M_0), \lambda_{M_0})$. Using the map $\cH^{\CAT}(M_0,M_1, h)\to \hAut (M_0, \del M_0)$ described above, we obtain a map \[\cH^{\CAT}(M_0, M_1,h)\to \cV(H_2(M_0), \lambda_{M_0}), \qquad W \mapsto \Delta_{F\circ \varphi_{W}}.\]

For our second invariant, let $M$ be a spin $4$-manifold and let $f$ be a homotopy equivalence restricting to the identity on $\del M$. Given a fixed spin structure $\fs$, there exists a free and transitive action of $H^1(M, \del M;\ZZ/2)$ on the set of spin structures restricting to $\fs\vert_{\del M}$ on $\del M$. Because $f$ fixes the boundary, $f^*(\fs)$ is a spin structure such that $f^*(\fs)\vert_{\del M}=\fs\vert_{\del M}$. Hence, we define $\Theta(f)$ to be the group element of $H^1(M,\del M;\ZZ/2)$ corresponding to the action taking $\fs$ to $f^*(\fs)$. Similarly to for variations, we prove that this provides a map \[\cH^{\CAT}(M_0, M_1,h)\to H^1(M_0,\del M_0;\ZZ/2), \qquad W \mapsto \Theta(F\circ \varphi_{W}).\] 

With this in hand, we can now state our main theorem. 
\begin{thm}\label[thm]{thm}
    Let $M_0, M_1$ be homeomorphic compact simply-connected smooth (or topological) $4$-manifolds. Fix a homeomorphism $F\maps M_1 \to M_0$ with $F\vert_{\del M_1}=h^{-1}$ for some homeomorphism $h\maps \del M_0 \to \del M_1$. 
    \begin{enumerate}
        \item When $M_0, M_1$ are not spin, there exists a bijection \[\cH^{\Diff}(M_0,M_1,h)\to \cV(H_2(M_0), \lambda_{M_0})\] given by $W \mapsto \Delta_{F\circ \varphi_W}$.
        \item When $M_0, M_1$ are spin, there exists a bijection \[\cH^{\Diff}(M_0,M_1,h)\to H^1(M_0, \del M_0;\ZZ/2) \times \cV(H_2(M_0), \lambda_{M_0})\] given by $W \mapsto (\Theta(F \circ \varphi_W),\Delta_{F\circ \varphi_W})$.
    \end{enumerate}
    In the case where $M_0= M_1$ and $h=\id$, we take $F=\id$ and these bijections are group isomorphisms. The same statement holds with $\cH^{\Diff}$ replaced by $\cH^{\Top}$. 
\end{thm}

\begin{rmk}
    When $M_0$ and $M_1$ are closed, this result recovers the result of \cite{Kreck_hcob}. To see this, we first note that for $M_0$ closed, $H^1(M_0, \del M_0;\ZZ/2)\cong H^1(M_0;\ZZ/2)$, which is trivial. Hence, in both cases, $h$-cobordisms are in bijection with Poincar\'e variations. Next, denoting $j_*\maps H_2(M_0)\to H_2(M_0, \del M_0)$ to be the map induced by the inclusion, \cite[Lemma 3.4]{Saeki} shows $\id-\Delta\circ j_*$ is an isometry and defines a map
    \[\Xi\maps\cV(M_0,\del M_0)\to \Aut(H_2(M_0),\lambda_{M_0}), \qquad \Delta\mapsto \id-\Delta\circ j_*,\]
    where $\Aut(H_2(M_0),\lambda_{M_0})$ is the group of isometries of the intersection form on $M_0$. In \cite[Theorem 7.13]{OrsonPowell} (see also \cite[Propositions 4.2, 4.3]{Saeki}), they show that this fits into an exact sequence
    \[0\to \wedge^2 H^1(\del M_0)^*\to \cV(M_0,\del M_0)\xrightarrow{\Xi}\Aut_\del^{\text{fix}}(H_2(M_0), \lambda_{M_0})\to 0,\]
    where $\Aut_\del^{\text{fix}}(H_2(M_0), \lambda_{M_0})$ is a subgroup of $\Aut(H_2(M_0), \lambda_{M_0})$. Here, we have $\wedge^2 H^1(\del M_0)^*=0$ and $\Aut_\del^{\text{fix}}(H_2(M_0), \lambda_{M_0})=\Aut(H_2(M_0), \lambda_{M_0})$ because $\del M_0=\emptyset$, so $\Xi$ is an isomorphism and we get a bijection $\cH^{\CAT}(M_0, M_1,h)\to \Aut(H_2(M_0),\lambda_{M_0})$. In particular, $\Delta_f\circ j_*=\id-(F\circ \varphi_W)_*$ so this bijection is given by $W\mapsto (F\circ \varphi_W)_*$. 
\end{rmk}

\begin{rmk}
    Given $M_0, M_1, h$, Boyer provides criteria in \cite[Theorem 0.7, Proposition 0.8]{Boyer} for the existence of a homeomorphism $F$ with $F\vert_{\del M_1}=h^{-1}$. We remark that this is necessary and sufficient for $\cH^{\Top}(M_0,M_1,h)$ to be nonempty. As mentioned above, given an $h$-cobordism in $\cH^{\Top}(M_0,M_1,h)$, the $h$-cobordism theorem provides such an $F$. For sufficiency, given such a homeomorphism $F$, the product cobordism $M_0 \times [0,1]$, lies in $\cH^{\Top}(M_0,M_1,h)$ if we take the inclusion of $M_0$ into $M_0 \times \{0\}$ to be the identity and the inclusion of $M_0$ into $M_0\times \{1\}$ to be $F$. It is also necessary and sufficient for $\cH^{\Diff}$ to be nonempty: the same argument shows necessity while for sufficiency, we can use the construction in \cref{sm_surj}. 
\end{rmk}

The main ingredient in our proof is Kreck's modified surgery program (\cite{SAD}), which we use to prove the injectivity portion of the theorem. We remark that our result recovers \cite[Theorem 3.1]{OrsonPowell}, which is the injectivity portion of their computation of the mapping class group for simply-connected $4$-manifolds with boundary; we spell this out in \cref{pseudoisotopies}. In their proof, they also used modified surgery but in a different flavor. While they work with the normal $2$-type, we work with the normal $1$-type, which reduces the amount of algebraic topology needed at the cost of introducing additional work to analyze the modified surgery obstruction. 

To be specific, in \cite[Proposition 6.2]{OrsonPowell}, they have a $4$-manifold $M$ and need to extend a particular normal $2$-smoothing on $\del M \times I$ to one on $M \times I$, where the smoothing was defined using a homeomorphism $G\maps M\to M$. Via an obstruction-theoretic argument, they show that $\Theta(G)$ appears as the primary obstruction to doing so and $\Delta_G$ as the secondary. An important element of their proof was using the product cell structure on $M\times I$ to describe the obstruction cocycle concretely. 

In our setting, we don't have a nice cell structure on any given $h$-cobordism so we instead utilize the $1$-type, for which the analogous step is constructing a normal $1$-smoothing on the union $W \cup W'$ where $W,W'\in \cH^{\CAT}(M_0,M_1,h)$. This amounts to finding a spin structure on $W \cup W'$, which is controlled by $\Theta (F \circ \varphi_W)$. The cost of this is that, unlike in \cite{OrsonPowell}, our surgery obstruction does not live inside a trivial group (specifically $L_6(\ZZ,S=\ZZ)$); instead, it lives in the monoid $\ell_6(\ZZ)$ defined in \cite{SAD}, so, we have to work to show that the obstruction is elementary. Essential to this is a computation of the homology groups of $W \cup W'$, which is the main place where $\Delta_{F \circ \varphi_{W}}$ appears (we also use this computation to show that in the nonspin case, $W \cup W'$ is oriented nullbordant). 

\subsection*{Organization}
We start by giving precise definitions of $h$-cobordisms for manifolds with boundary and their invariants in \cref{invts}. Then, we proceed to briefly review classical and modified surgery theory in \cref{surgery}. Additionally, in \cref{modified}, we simplify the analysis of the surgery obstruction in particular cases to a homology computation. Afterwards, we proceed to prove the surjectivity and injectivity of \cref{thm} in \cref{existence} and \cref{uniqueness} respectively. We compile these results to prove the main theorem in \cref{proof} and 
conclude by connecting the content herein to pseudoisotopies in \cref{pseudoisotopies}. 

\subsection*{Conventions}
Unless stated otherwise, all manifolds will be assumed to be compact, connected, and orientable. For $f$ a homotopy equivalence, we will denote by $f^{-1}$ a homotopy inverse of $f$. Also, when we want to work in either the smooth or topological category, we utilize the notation $\CAT$ to represent $\Diff$ or $\Top$. 

\subsection*{Acknowledgments}
I would like to thank my advisor, Anthony Conway, for helpful advice, conversations, and support. I would also like to thank Audrick Pyronneau and Daniel Galvin for helpful discussions. 

\subsection*{AI Statement}
No AI tools were used in the preparation of this paper. 

\section{Invariants of $h$-corbordisms rel. boundary}\label[section]{invts}
\subsection{$h$-cobordisms}\label[section]{hcob}
We briefly recall a precise definition of $h$-cobordism here because the level of detail will be useful to us throughout (see also \cite[Chapter 1]{Luck}, \cite[Section 4.1]{OrsonPowell}). 
\begin{defn}
    Let $M_0,M_1$ be oriented $n$-dimensional smooth manifolds with $\del M_0 \cong \del M_1$ via an orientation preserving diffeomorphism $h\maps \del M_0 \to \del M_1$. An $(n+1)$-dimensional \emph{$h$-cobordism rel. h} from $M_0$ to $M_1$ is a tuple $(W,M_0,f_0,M_1,f_1,g)$ where 
    \begin{enumerate}[i)]
        \item W is an $(n+1)$-dimensional smooth manifold with corners;
        \item $\del W = \del_- W \cup \del_V W \cup \del_+ W$, where the interiors of $\del_- W, \del_V W,$ and $ \del_+ W$ form the codimension $1$ stratum of $W$ and $\del\del_-W, \del\del_+W$ form the codimension $2$ corners;
        \item $f_0\maps (M_0,\del M_0) \to (\del_- W,\del\del_-W)$, $g\maps \del M_0 \times I \to \del_V W$, and $f_1\maps (\overline{M_1}, \overline {\del M_1}) \to (\del_+ W, \del\del_+W)$ are orientation preserving diffeomorphisms, with $\del W$ possessing the boundary orientation from $W$;
        \item The inclusion maps $i_\pm \maps (\del_k W, \del \del_\pm W) \to (W, \del_V W)$ are homotopy equivalences of pairs;
        \item The following diagrams commute, where $i_V\maps \del_V W\to W$ is the inclusion: 
        \[\begin{tikzcd}
\del M_0 \arrow[d, "f_0"'] \arrow[r, "g"] & \del_V W \arrow[d, "i_V"] &  & \del M_1 \arrow[r, "g\circ h^{-1}"] \arrow[d, "f_1"'] & \del_V W \arrow[d, "i_V"] \\
\del_- W \arrow[r, "i_-"]                 & W                         &  & \del_+ W \arrow[r, "i_+"]                             & W                        
\end{tikzcd}\]

    \end{enumerate}
The same definition applies in the topological category, except we don't need the decomposition of $W$ as a manifold with corners and instead view it as a manifold with boundary. 
\end{defn}
In the introduction, we said that an $h$-cobordism rel. $h$ has a boundary that decomposes as $M_0\cup_{\id} (\del M_0 \times [0,1])\cup_h M_1$; condition v) guarantees that $f_0\cup g \cup f_1$ gives the desired decomposition. 
\begin{rmk}
    In this definition, we keep track of the maps $f_0, f_1,g$ because it will be necessary to keep track of how their boundaries are identified with $M_0 \cup \del M_0\times [0,1] \cup M_1$, e.g. for gluing purposes. 
\end{rmk}
\begin{rmk}\label[rmk]{replace_g}
    For notational reasons, we usually refer to the cobordism $(W,M_0,f_0,M_1,f_1,g)$ just as the $(n+1)$-manifold $W$, and for a map $g\maps \del M_0 \times [0,1] \to W$, we denote the restriction to $\del M_0 \times \{i\}$ as $g_i$. When appropriate, we also abuse notation and omit some of the inclusions, writing $i_V\circ g=g$ as well as $i_\pm^{-1}\circ i_V \circ g_i = i_V \circ g_i= g_i$. 
\end{rmk}

\begin{rmk}
    Given $(W,M_0,f_0,M_1,f_1,g)$ an $h$-cobordism rel. $h$ and a $\CAT$-isomorphism $g'$ such that $f_0^{-1} \circ g'_0=\id$, we can also form an $h$-cobordism $(W,M_0,f_0,M_1,f_1,g')$. This gives us an $h$-cobordism rel. $h'$, where $h'$ is a $\CAT$-isomorphism pseudoisotopic to $h$ via the map 
    \[(h\times \id)\circ g^{-1}\circ g'\maps \del M_0 \times [0,1] \to \del M_1 \times [0,1].\]
    This restricts to $h$ on $\del M_1 \times\{0\}$ and on $\del M_1 \times \{1\}$, it restricts to
    \[h\circ g_1^{-1}\circ g_1'=(g_1\circ h^{-1})\circ g_1'=f_1^{-1}\circ g_1'=h',\]
    so it does indeed give the desired pseudoisotopy. 
\end{rmk}

\begin{comment}
\begin{rmk}
    If we only require $\del_V W$ to be homotopy equivalent to $\del M_0 \times I$ then we call $W$ a relative $h$-cobordism, motivated by the conventions in eg \cite{ranickiexactseq} and \cite{structuresets}. 
\end{rmk}
\end{comment}

Similarly, we also recall a precise definition of what it means for two cobordisms rel. $h$ to be $\CAT$-isomorphic rel. boundary. Define $\psi_W:=f_0 \cup g \cup f_1$. Then, we have the following definition. 

\begin{defn}
    Let $(W,M_0,f_0,M_1,f_1,g)$ and $(W',M_0,f_0',M_1,f_1',g')$ be $h$-cobordisms rel. $h$. Then, $W$ and $W'$ are \emph{$\CAT$-isomorphic rel. boundary} if there exists a $\CAT$-isomorphism $\Phi\maps W \to W'$ such that the following diagram commutes:
    \[
\begin{tikzcd}
{M_0 \cup_{\id} \del M_0 \times [0,1] \cup_{h} M_1} \arrow[r, "\psi_W"] \arrow[d, "\id"] & {\del_- W \cup_{\id} \del_V W \cup_{\id} \del_+ W} \arrow[d, "\Phi\vert_{\del W}"] \\
{M_0 \cup_{\id} \del M_0 \times [0,1] \cup_{h} M_1} \arrow[r, "\psi_{W'}"']              & {\del_- W' \cup_{\id} \del_V W' \cup_{\id} \del_+ W'} .                             
\end{tikzcd}\]
\end{defn}

\begin{defn}\label[defn]{stack}
    We denote the set of $\CAT$ $h$-cobordisms rel. $h$ from $M_0$ to $M_1$ to be $\cH^{\CAT}(M_0, M_1,h)$.
\end{defn}
We also introduce some notation for concatenating $h$-cobordisms. 
\begin{defn}\label[defn]{mult}
    Let $M_0,M_1,M_2$ be $n$-dimensional $\CAT$-manifolds with $\del M_0 \cong \del M_1 \cong \del M_2$. Define an operation \[*\maps \cH^{\CAT}(M_1, M_2,h') \times \cH^{\CAT}(M_0, M_1,h)\to \cH^{\CAT}(M_0, M_2,h'\circ h)\] in the following manner:
    
    Let $(W,M_0,f_0,M_1,f_1,g)$ and $(W',M_1,f_1',M_2,f_2,g')$ be $h$-cobordisms rel. $h, h'$ respectively and define $Y=W\cup_{f_0'\circ f_1^{-1}} W'$, which intuitively speaking, is $W'$ stacked on top of $W$. Observe that \[\del Y = \del_-W \cup_{\id} \del_V W \cup_{f_1'\circ f_1^{-1}} \del_V W' \cup_{\id} \del_+W'.\] 
    Define the map $\tilde g \maps \del M_0\times [0,1] \to \del_V W \cup_{f_1'\circ f_1^{-1}} \del_V W'$ via 
    \[\tilde g(x,t)=
    \begin{cases}
        g(x,2t)&x\in[0,1/2]\\
        g'((g_0'^{-1}\circ f_1'\circ f_1^{-1}\circ g_1) (x),2t-1) &x\in[1/2,1].
    \end{cases}
    \]
    This is well-defined because $g_0'=g'\vert_{\del M_0\times \{0\}}$ so at $t=1/2$, we have that $g'$ and $g_0'^{-1}$ cancel. So, denote $\del_V Y := \del_V W \cup_{f_1'\circ f_1^{-1}} \del_V W'$ and define $\tilde g_i$ like before. Observe that $\tilde g(x,0)=f_0(x)$ and
    \[f_2^{-1}\circ \tilde g_1=f_2^{-1}\circ g'_1\circ g_0'^{-1}\circ f_1'\circ f_1^{-1}\circ g_1=h' \circ h,\]
    noting that $f_1'\vert_{\del M_0}=g_0'$ so $g_0'^{-1}\circ f_1'\vert_{\del M_0} = \id$. Hence, we can define \[W'*W:=(Y, M_0, f_0, M_1, f_2, M_2, \tilde g),\] which is an $h$-cobordism rel. $h'\circ h$. 
\end{defn}

Let $\HE(M_0,M_1,h)$ be the set of orientation preserving homotopy equivalences $M_0 \to M_1$ restricting to $h$ on the boundary, modulo homotopy rel. $\del M_0$. Note that for $h$ an orientation preserving homeomorphism, any homotopy equivalence restricting to $h$ on the boundary is automatically orientation preserving so the orientation preserving condition is only relevant when $M_0, M_1$ are closed.

Given an $h$-cobordism $(W,M_0,f_0,M_1,f_1,g)$ between closed manifolds, there is a standard way to obtain a homotopy equivalence, namely via the map \[f_1^{-1}\circ i_+^{-1}\circ i_- \circ f_0 \maps M_0 \to M_1,\] 
where $i_+^{-1}$ is any choice of homotopy inverse. In the case of an $h$-cobordism rel $h$, we would like the same construction to produce a map $\cH^{\CAT}(M_0,M_1,h) \to \HE(M_0,M_1,h)$. In order to do so, we need to construct a particular $i_+^{-1}$ because a priori, we don't know if $f_1^{-1}\circ i_+^{-1}\circ i_- \circ f_0$ restricts to $h$ on $\del M_0$. We start this process by identifying the condition we need $i_+^{-1}$ to satisfy. 

\begin{lemma}\label[lemma]{phi_cond}
    Let $i_+^{-1}$ be a deformation retract providing a homotopy inverse for $i_+$. The map \[f_1^{-1}\circ i_+^{-1}\circ i_- \circ f_0 \maps M_0 \to M_1,\] restricts to $h$ on $\del M_0$ if and only if $(i_+^{-1}\circ g) (x,0)=g(x,1)$ for all $x\in \del M_0 \times [0,1]$. 
\end{lemma}
\begin{proof}
Recalling that we usually suppress $i_-$ from our notation and that $f_0\vert_{\del M_0}=g_0$, we get \[f_1^{-1} \circ i_+^{-1}\circ i_- \circ f_0\vert_{\del M_0} =f_1^{-1}\circ i_+^{-1}\circ g_0.\] By definition, $h=f_1^{-1} \circ g_1$, so $f_1^{-1}\circ i_+^{-1}\circ i_- \circ f_0\vert_{\del M_0}=h$ if and only if $f_1^{-1}\circ i_+^{-1}\circ g_0=f_1^{-1} \circ g_1$, which is equivalent to the statement $i_+^{-1} \circ g_0=g_1$. Recalling that $g_i(x):=g(x,i)$, we see that $i_+^{-1}\circ g_0=g_1$ if and only if $(i_+^{-1}\circ g) (x,0)=g(x,1)$ for all $x\in \del M_0 \times [0,1]$. 
\end{proof}

In light of the above, we will construct a deformation retraction $H$ from $W$ to $\del_+ W$ such that $H(-,1)$ satisfies the conditions of \cref{phi_cond}. In fact, we construct it so that it satisfies something stronger, namely that elements of $\del_V W$ remain in $\del_V W$ throughout the homotopy; this condition will be needed later in the proof of \cref{def_unique}. 

\begin{lemma}\label[lemma]{deformation_construction}
There exists a deformation retraction $H\maps W \times I \to W$ from $W$ to $\del_+ W$ such that for all $(y,z)\in \del M_0 \times [0,1]$ and $t\in[0,1]$, there exists $z' \in [0,1]$ such that $H(g(y,z),t))=g(y,z')$
\end{lemma}
\begin{proof}
We will construct the desired deformation retraction by first defining the deformation retraction on $\del_V W \cup \del_+ W$, then extending it to a homotopy on $W$, which we will then concatenate with an arbitrary deformation retraction to get the desired map. 

For elements in $\del_V W$, we will use the identification with $\del M_0 \times [0,1]$ to write them as $g(y,z)$, with $(y,z)\in \del M_0 \times [0,1]$. Consider the homotopy $H_0 \maps \del_V W \cup \del_+ W \to \del _+ W$ given by 
\[H_0(g(y,z),t))=g(y,t+(1-t)z),\] 
on $\del_V W$ and the constant homotopy on $\del_+ W$. 

Note that $(W, \del_V W \cup \del_+ W)$ satisfies the homotopy extension property because $\del_V W  \cup \del_+ W$ is locally flatly embedded and hence, has a collar neighborhood by \cite[Theorem 2]{Brown}. So, we can extend it to a homotopy $H_1\maps W \times[0,1] \to \del_V W$ such that $H_1\vert_{W\times\{0\}}=\id$ and $H_1\vert_{\del_{V} W \cup \del_{+} W \times [0,1]}=H_0$.

For the same reason, $(W, \del_+ W)$ satisfies the homotopy extension property. So, $i_+$ has a homotopy inverse that is a deformation retract, e.g. via \cite[Corollary 0.20]{Hatcher}. Thus, let $H_2 \maps W \times [0,1] \to W$ from $W$ to $\del_+ W$ be a deformation retraction from $W$ to $\del_+ W$. Define
\[H(x,t)=\begin{cases}
    H_1(x,2t) & t\in [0,1/2] \\
    H_2(H_1(x,1),2t-1) &t \in [1/2,1].
\end{cases}\]
To see that this is a deformation retraction, observe that
\begin{enumerate}[i)]
    \item $H\vert_{W\times\{0\}}=H_1\vert_{W\times\{0\}}=\id$,
    \item $H(x,1)\in \del_+W$ for all $x\in W$ because $H_2$ was a deformation retraction, and
    \item $H(x,t)=x$ for all $x\in \del_+ W$ because $H_1$ restricts to the constant homotopy on $\del_+W$ and $H_2$ was a deformation retraction. 
\end{enumerate}
In addition, for $(y,z)\in \del M_0 \times [0,1]$, 
\[H(g(y,z),t)=\begin{cases}
    g(y,2t+(1-2t)z) &t\in [0,1/2]\\
    g(y,1) &t\in [1/2,1],
\end{cases}\]
so in both cases, the statement of the theorem is satisfied. 
\end{proof}

A priori, there are many different ways of finding a deformation retract satisfying the conditions of \cref{phi_cond}. Despite this, we show that any two choices are in fact homotopy equivalent rel. $\del \del_- W$.

\begin{lemma}\label[lemma]{def_unique}
    Any two deformation retractions $H_0, H_1\maps W \times [0,1] \to W$ from $W$ to $\del_+ W$ such that $H_i(g(y,0),1)=g(y,1)$ for all $y\in \del M_0$ are homotopic rel. $\del \del_- W \times \{1\}$. Furthermore, we can choose this homotopy $H_s$ such that for all $s\in[0,1]$, $H_s$ is a deformation retraction onto $\del_+ W$. 
\end{lemma}
\begin{proof}
    It suffices to prove this under the assumption that $H_1$ is a deformation retraction from $W$ to $\del_+ W$ such that for $(y,z)\in \del M_0 \times [0,1]$, and all $t\in[0,1]$, we have $H_1(g(y,z),t))=g(y,z')$ for some $z' \in [0,1]$. This exists via \cref{deformation_construction}. Define a homotopy from $H_0$ to $H_1$ via
    \[H_s(x,t)=\begin{cases}
        H_0(H_1(x,2st),t) & s\in [0,1/2]\\
        H_0(H_1(x,t),2t(1-s)) & s\in[1/2,1].
    \end{cases}\]
    For $x \in \del \del_- W $, we can write $x=g(y,0)$ for some $y \in \del M_0$. Then, for $s\in[0,1/2]$, there exists $z' \in [0,1]$ such that
    \[H_s(x,1)=H_0(H_1(g(y,0),2s),1)=H_0(g(y,z'),1)=g(y,1).\]
    Meanwhile, for $s\in[1/2,1]$, note that $H_1(g(y,0),1)=g(y,1)$. This is because $H_1(g(y,0),1))=g(y,z')$ for some $z' \in [0,1]$, and $H_1$ is a deformation retraction onto $\del_+ W$ so $g(y,z')\subset \del \del_+ W$, meaning that in fact, $z'=1$. Additionally, because $H_0$ is a deformation retraction, it fixes $g(y,1)$, so we obtain
    \[H_s(x,1)=H_0(H_1(g(y,0),1),2-2s)=H_0(g(y,1),2-2s)=g(y,1).\]
    Hence, $H_s$ is a homotopy rel. $\del \del_- W \times \{1\}$. Hence, $H_0$ and $H_1$ are homotopic rel. $\del \del_- W \times \{1\}$. A quick check of the conditions also shows that $H_s$ is a deformation retraction onto $\del_+W$ for all $s$. 
\end{proof}

Putting things together, we get the following proposition. 

\begin{prop}\label[prop]{varphi}
    Let $(W,M_0, f_0, M_1, f_1, g)$ be an $h$-cobordism rel. $h$. Then, there exists a well-defined map $\cH^{\CAT}(M_0, M_1, h) \to \HE(M_0, M_1, h)$ defined via \[W \mapsto \varphi_W:=f_1^{-1}\circ i_+^{-1} \circ i_- \circ f_0,\] where $i_+^{-1}$ is any deformation retract onto $\del_+W$ such that $(i_+^{-1}\circ g) (x,0)=g(x,1)$ for all $x\in \del M_0 \times [0,1]$. 
\end{prop}
\begin{proof}
    We start by noting that such an $i_+^{-1}$ exists via \cref{deformation_construction} and that $\varphi_W\vert_{\del M_0}=h$ via \cref{phi_cond} so $\varphi_W\in \HE(M_0, M_1, h)$. 

    To show that $\varphi_W$ is independent of choice of $i_+^{-1}$, let $j_0,j_1$ be deformation retracts onto $\del_+W$ such that $(j_i^{-1}\circ g) (x,0)=g(x,1)$ for all $x\in \del M_0 \times [0,1]$. Take deformation retractions $H_0, H_1$ such that $H_i(-,1)=j_i$. By \cref{def_unique}, $H_0$ and $H_1$ are homotopic rel. $\del \del_- W \times \{1\}$ via a homotopy $H_s$ such that $H_s$ is a deformation retract onto $\del_+W$ for all $s$. In particular, we obtain a homotopy
    \[\varphi_s:=f_1^{-1}\circ H_s\vert_{\del_-W\times\{1\}}\circ i_- \circ f_0\]
    from $f_1^{-1}\circ j_0\circ i_- \circ f_0$ to $f_1^{-1}\circ j_1\circ i_- \circ f_0$. Note that $\varphi_s$ is well defined because we chose $H_s$ to be a deformation retract for every $s$, guaranteeing that $H_s(\del_-W,1)\subset \del_+ W$. Additionally, for every $s$, $\varphi_s\in \HE(M_0,M_1,h)$ because $H_s$ was a homotopy rel. $\del\del_-W$. Therefore, the definition of $\varphi_W$ is unique up to homotopy rel. $\del M_0$. 
\end{proof}

For the remainder of this paper, we choose $i_+^{-1}$ to be any deformation retract onto $\del_+W$ such that $(i_+^{-1}\circ g) (x,0)=g(x,1)$ for all $x\in \del M_0 \times [0,1]$. The same construction provides a homotopy inverse $i_-^{-1}$ satisfying $g_0^{-1}\circ i_-^{-1}\circ g_1=\id$. So, we will also let $i_-^{-1}$ denote such a homotopy inverse.

When $M_0=M_1$ and $h=\id$, we denote $\HE(M_0, M_1,h)=\hAut(M_0, \del M_0)$ and give $\hAut$ a group structure via composition. We note here that elements of $\hAut(M_0, \del M_0)$ are invertible because any homotopy equivalence $f\in\hAut(M_0, \del M_0)$ has a homotopy inverse $f^{-1}$ restricting to the identity on $\del M_0$, and the homotopies $f\circ f^{-1}\simeq \id$ and $f^{-1}\circ f\simeq\id$ can be taken to be rel. $\del M_0$, e.g. via \cite[Section 6.5]{may}.

In this case, we obtain a group homomorphism $\cH^{\CAT}(M_0,M_0,\id)\to \hAut(M_0, \del M_0)$ as seen in the following lemma. 
\begin{lemma}\label[lemma]{hcob_composition}
    Let $(W,M_0,f_0,M_1,f_1,g)$ and $(W',M_1',f_1',M_1',f_2',g')$ be $h$-cobordisms. Then, \[\varphi_{W'*W}\simeq\varphi_W'\circ \varphi_{W}.\]
\end{lemma}
\begin{proof}
    Denote $Y=W'* W$ and let $i_\pm^Y$ be the inclusions $\del_\pm Y \to Y$ and likewise for $W,W'$. Similarly, let $(i_+^Y)^{-1},(i_+^W)^{-1},(i_+^{W'})^{-1}$ be homotopy inverses as in \cref{varphi}. Recalling that we chose these to be deformation retracts, let $H_W\maps W \times I \to W$ and $H_{W'}\maps W' \times I \to W'$ be the corresponding deformation retractions so that $(i_+^W)^{-1}(x)=H_W(x,1)$ and likewise for $W'$.
    Define $H \maps Y\times I \to Y$ via
    \[
    H(x,t)=\begin{cases}
        H_W(x,2t) & x\in W, t\in [0,1/2] \\
        x &x\in W', t\in [0,1/2]\\
        H_{W'} (i_-^{W'}\circ f_1' \circ f_1^{-1}\circ (i_+^W)^{-1}(x), 2t-1) & x\in W, t\in [1/2,1]\\
        H_{W'}(x,2t-1) &x \in W', t \in [1/2,1].
    \end{cases}
    \]
    To see this is well defined, note that for $x \in W \cap W'$, we have $H_W(x,2t)=x$ because $H_W$ was a deformation retraction and also, $i_-^{W'}\circ f_1' \circ f_1^{-1}\circ (i_1^W)^{-1}(x)=x$ because we glued $W$ to $W'$ via $f_1'\circ f_1^{-1}$. 

    So, we obtain a homotopy inverse $(i_+^Y)^{-1}:=H(-,1)$. A quick computation shows that $(i_+^Y)^{-1}$ satisfies $\tilde g^{-1}\circ (i_+^Y)^{-1}\circ \tilde g (x,0)=(x,1)$ for $x \in \del M_0$. Hence, by \cref{varphi}, $\varphi_Y\simeq f_2^{-1} \circ (i_+^Y)^{-1} \circ i_-^Y \circ f_0$. Thus,
    \begin{align*}
        \varphi_Y&\simeq f_2^{-1} \circ (i_+^Y)^{-1} \circ i_-^Y\circ f_0\\&=(f_1')^{-1}\circ(i_W^+)^{-1}\circ i_-^{W'}\circ f_1' \circ f_1^{-1}\circ (i_+^W)^{-1}\circ i_-^Y\circ f_0\\
        &=\varphi_W' \circ f_1^{-1}\circ (i_+^W)^{-1}\circ i_-^W\circ f_0\\
        &=\varphi_W'\circ \varphi_W,
    \end{align*}
    as desired. 
\end{proof}
\subsection{Variations}\label[section]{variations}
In this section, we define the first of the two invariants associated to an $h$-cobordism. 
\begin{defn}
    Given a continuous map of pairs $f\maps (M,N) \to (M,N)$ restricting to the identity on $N$, we can define the \textit{variation} associated to $f$ to be the map $\Delta_f \maps H_i(M, N) \to H_i(M)$ defined via $[\Sigma] \mapsto [\Sigma-f(\Sigma)]$. This is well defined (see e.g. \cite[Section 2.2]{OrsonPowell}). In the case that $M$ is a $4$-manifold, we can further define the set $\cV(H_2(M), \lambda_M)$ of \textit{Poincar\'e variations} to be homomorphisms $\Delta \maps H_2(M, \del M) \to H_2(M)$ such that 
\[\Delta+\Delta^!=\Delta \circ j_* \circ \Delta^!,\]
where $j_* \maps H_2(M) \to H_2(M, \del M)$ is the inclusion and $\Delta^!:=\PD\circ \ev^{-1} \circ \Delta^* \circ \ev\circ \PD^{-1}$, where $\ev$ and $\PD$ denote the evaluation map and Poincar\'e duality respectively. 
\end{defn}
Saeki proved that Poincar\'e variations form a group via the following. 
\begin{lemma}[{\cite[Lemma 3.5]{Saeki}}]
The set $\cV(M, \lambda_M)$ forms a group under the operation
\[\Delta_1 * \Delta_2=\Delta_1+(\id-\Delta_1\circ j_*)\circ \Delta_2.\]
\end{lemma}
In the case where $M$ is smooth and $f$ is a stable diffeomorphism, Saeki proves that $\Delta_f$ is a Poincar\'e variation (\cite[Lemma 3.2(4)]{Saeki}), and Orson-Powell use that the same holds when $f$ is a homeomorphism. We check that the same holds when $f$ is a homotopy equivalence. To do so, we require the following lemma from \cite{OrsonPowell}, which will be of further use later as well. 
\begin{lemma}[{\cite[Lemma 2.14]{OrsonPowell}}]\label[lemma]{htpy}
 Let $(A,B)$ be a pair of spaces and let $f,g \maps (A,B)\to (A,B)$ be continuous maps each restricting to the identity on $B$. If $f\simeq g$ via a homotopy restricting to the identity on $B$,  then $\Delta_{f}=\Delta_{g}$.
\end{lemma}

With this in hand, we prove the following:
\begin{lemma}\label[lemma]{htpy_gives_var}
    Let $f\maps M \to M$ be an orientation preserving homotopy equivalence restricting to the identity on the boundary. Then, the variation $\Delta_f$ is Poincar\'e.
\end{lemma}
\begin{proof}
    As in \cite[Lemma 3.2(4)]{Saeki}, the problem of checking whether the variation $\Delta_f$ is Poincar\'e reduces to checking for any $\alpha, \beta\in H_2(M, \del M)$ that $-Q_M^\del((f_*^{-1}\circ\Delta_f)(\alpha), \beta)=Q_M^\del(\Delta_f (\beta), \alpha)$, noting that \cite[Lemma 3.2(2)]{Saeki}, holds without modification. Here, $Q_M^\del$ denotes the relative intersection form $H_2(M)\times H_2(M, \del M) \to \ZZ$. We also use $Q_M$ to denote the intersection form $H_2(M)\times H_2(M) \to \ZZ$. 
    Let $a,b$ be surfaces with boundary representing $\alpha, \beta$ respectively and let $f^{-1}$ be any homotopy inverse to $f$ restricting to the identity on the boundary. Then, 
    \begin{align*}
        -Q_M^\del((f_*^{-1}\circ\Delta_f)(\alpha), \beta)
        &=-Q_M^\del([f^{-1}(a)-f^{-1}\circ f(a)], \beta)\\
        &=-Q_M^\del([-(a-f^{-1}(a))+a-f^{-1}\circ f(a)], \beta)\\
        &=Q_M^\del(\Delta_{f^{-1}}(\alpha)-\Delta_{f^{-1}\circ f}(\alpha),\beta).
    \end{align*}
    Noting that $f^{-1}\circ f\simeq \id$ rel. boundary, \cref{htpy} provides $\Delta_{f^{-1}\circ f}=\Delta_{\id}$, which is the zero map. Hence, 
    \[Q_M^\del(\Delta_{f^{-1}}(\alpha)-\Delta_{f^{-1}\circ f}(\alpha),\beta)=Q_M^\del(\Delta_{f^{-1}}(\alpha),\beta).\]
    
    It will be convenient to work with absolute classes instead of relative ones because both versions of the intersection form require at least one of the classes to be absolute. To do so, we define an auxiliary manifold $\hat M$ by attaching $2$-handles along each connected component of $\del a$ and $\del b$. Then, let $\hat a$ be the surface obtained by taking the union of $a$ with the core of the $2$-handles attached to $\del a$ and $\hat b$ likewise, with $\hat \alpha$ and $\hat \beta$ being the respective absolute homology classes; by a slight abuse of notation, we also use the same notation when viewing them as relative homology classes. Also define $\hat f\maps \hat M \to \hat M$ by extending $f$ over the $2$-handles via the identity and $\hat f^{-1}$ similarly; this is well defined because $f$ and $f^{-1}$ are the identity on $\del M$. Note that $\hat f^{-1}$ still provides a homotopy inverse for $f$ because we can extend any homotopy across the $2$-handles using the identity. 

    Let $c$ be a surface representing $\Delta_{f^{-1}}(\alpha)$ in $H_2(M)$ and intersecting $b$ transversely. Then, $c$ intersects $\hat b$ transversely in $\hat M$ because all intersection points lie inside $M$, within which $\hat b$ and $b$ agree. In particular, the signed count of the intersections between $c$ and $b$ is the same as the signed count between $c$ and $\hat b$. Hence, we obtain
    \[Q_M^\del(\Delta_{f^{-1}}(\alpha),\beta)=Q_{\hat M}(\Delta_{f^{-1}}(\alpha),\hat \beta).\]
    
    Recall that $Q_M^\del$ is defined via $\PD^{-1}\circ \ev \maps H_2(M,\del M)\to H_2(M)^*$ while $Q_M$ is defined similarly by $i_* \circ \PD^{-1}\circ \ev \maps H_2(M)\to H_2(M)^*$, where $i_*$ is the map $H_2(M)\to H_2(M, \del M)$ induced by the inclusion. Therefore,
    \[Q_{\hat M}(\Delta_{f^{-1}}(\alpha),\hat \beta)=Q_{\hat M}^\del(i_*\circ\Delta_{f^{-1}}(\alpha),\hat \beta).\]
    
    Letting $d$ be the union of the cores of the $2$-handles attached to $\del a$, we obtain
    \[i_*\circ\Delta_{f^{-1}}(\alpha)=[a-f^{-1}(a)]=[a+d-f^{-1}(a)-d]=[\hat a-\hat f^{-1}(\hat a)]=[\hat a]-[\hat f^{-1}(\hat a)]=\hat \alpha-\hat f^{-1}(\hat \alpha).\]
    Here, we note that we had to switch from $a$ to $\hat a$ because $a$ is not a relative cycle in $C_2(\hat M,\del \hat M)$ (because $\del a$ is contained in the interior of $\hat M$) but $\hat a$ is an absolute cycle, so in particular, it is a relative cycle in $C_2(\hat M,\del \hat M)$. So, we have
\begin{align*}
        Q_{\hat M}^\del(i_*\circ\Delta_{f^{-1}}(\alpha),\hat \beta)
        =Q_{\hat M}^\del(\hat\alpha-\hat f^{-1}_*(\hat\alpha),\hat \beta)
        =Q_{\hat M}^\del(\hat\alpha,\hat \beta)-Q_{\hat M}^\del (\hat f^{-1}_*(\hat\alpha),\hat \beta).
    \end{align*}
Because $\hat f$ is also an orientation-preserving homotopy equivalence, and hence an isometry, we obtain
    \begin{align*}
        Q_{\hat M}^\del(\hat \alpha,\hat \beta)-Q_{\hat M}^\del (\hat f^{-1}_*(\hat \alpha),\hat \beta)
        &=Q_{\hat M}^\del(\hat \alpha,\hat \beta)-Q_{\hat M}^\del (\hat \alpha,\hat f_*(\hat \beta))\\
        &=Q_{\hat M}^\del(\hat \alpha,\hat \beta -\hat f_*(\hat \beta))\\
        &=Q_{\hat M}(\hat \alpha,\hat \beta -\hat f_*(\hat \beta))\\
        &=Q_{\hat M}(\hat \beta -\hat f_*(\hat \beta),\hat \alpha),
    \end{align*}
    where when switching from $Q_{\hat M}^\partial$ to $Q_{\hat M}$, we view $\hat \alpha$ and $\hat \beta$ as absolute instead of relative classes. By similar arguments as before, we conclude
    \[Q_{\hat M}(\hat \beta -\hat f_*(\hat \beta),\hat\alpha)=Q_{\hat M}(\Delta_f(\beta),\hat\alpha)=Q_M^\del(\Delta_f(\beta),\alpha),\]
    as required.
 \end{proof}

In the case of homeomorphisms, \cite[Corollary 2.15]{OrsonPowell} shows that the map \[\tilde\pi_0(\Homeo^+(M,\del M))\to \cV(H_2(M),\lambda_M), \qquad f\mapsto \Delta_f,\] 
is a group homomorphism, where $\tilde \pi_0 (\Homeo^+(M,\del M))$ is the group of orientation preserving homeomorphisms of $M$ restricting to the identity on the boundary modulo pseudoisotopy. The same proof shows that the map $\hAut(M,\del M)\to \cV(H_2(M),\lambda_M)$ defined in the same manner is also a group homomorphism. For completeness, we perform the proof here. 

\begin{lemma}\label[lemma]{var_hom}
    The map $\hAut(M,\del M)\to \cV (H_2(M),\lambda_M)$ given by $f\mapsto \Delta_f$ a group homomorphism.
\end{lemma}
\begin{proof}
    By \cref{htpy_gives_var,htpy}, this is well defined. Furthermore, \cite[Lemma 3.2(5)]{Saeki} goes through without change so given $f,g \in \hAut(M,\del M)$,
    \[\Delta_{f\circ g}=\Delta_f+\Lambda_f\circ \Delta_g=\Delta_f+(\id-\Delta_f \circ j_*)\circ \Delta_g=\Delta_f *\Delta_g.\]
\end{proof}

We finish this subsection by connecting the content herein to $h$-cobordisms. 

\begin{lemma}\label[lemma]{hcob_var_hom}
    Let $M_0, M_1$ be homeomorphic simply-connected $\CAT$ $4$-manifolds. Fix a homeomorphism $F\maps M_1 \to M_0$ with $F\vert_{\del M_1}=h^{-1}$ for some $h$. Then, the map $\cH^{\CAT}(M_0, M_1, h)\to \cV(H_2(M_0), \lambda_{M_0})$ given by $W \mapsto \Delta_{F\circ \varphi_{W}}$ is a well-defined map. Furthermore, in the case where $M_0=M_1$ and $h=\id$, we choose $F=\id$ and this map is a group homomorphism. 
\end{lemma}
\begin{proof}
    In \cref{varphi}, we built a map $\cH^{\CAT}(M_0, M_1, h) \to \HE (M_0, M_1, h)$ mapping $W$ to $\varphi_W$. Postcomposing with $F$ provides a map $\HE(M_0, M_1, h)\to \hAut(M_0, \del M_0)$, and taking variations provides a map $\hAut(M_0, \del M_0)\to \cV(H_2(M_0), \lambda_{M_0})$. Hence, taking the composition of these three maps provides a well-defined map $\cH^{\CAT}(M_0, M_1,h) \to \cV(H_2(M_0), \lambda_{M_0})$ via $W\mapsto \Delta_{F\circ\varphi_W}$. 
    
    When $M_0 = M_1$ and $h=\id$, we pick $F=\id$ and this map is a group homomorphism via \cref{hcob_composition} and \cref{var_hom}. 
\end{proof}
\subsection{The spin invariant}\label[section]{spin_invt}
For the case where our $h$-cobordisms are spin, we have an additional invariant to measure whether they induce the same spin structures on their boundaries. We briefly review some background on relative spin structures here; for further details, see \cite[Section 2.1]{OrsonPowell}.

Recall that the stable normal bundle of a smooth manifold $M$ is classified by the stable Gauss map $\nu \maps M \to \BSO$, and that the universal covering map $\Spin \to \SO$ induces a principal fibration on classifying spaces $\gamma \maps \BSpin \to \BSO$. We say that $M$ is \textit{spin} if it admits a lift $\fs \maps M \to \BSpin$; we call such a lift a \textit{spin structure}. Two spin structures are \textit{isomorphic} if they are homotopic through spin structures. By obstruction theory, there is a free transitive action of $H^1(M;\ZZ/2)$ on the set of spin structures on $M$ up to isomorphism. In the topological category, there exists an analogous definition using topological stable normal bundles and replacing $\BSO, \BSpin$ with $\BTop, \BTopSpin$ respectively. 

When $M$ has boundary, a spin structure $\fs\maps M \to \BCATSpin$ restricts to provide a spin struture $\del \fs\maps \del M \to \BCATSpin$ on $\del M$. We say two spin structures $\fs_0$ and $\fs_1$ are \textit{isomorphic rel. boundary} if $\del \fs_0=\del \fs_1$ and $\fs_0 \simeq \fs_1$ through spin structures restricting to $\del \fs_0$ on the boundary. Like before, there is a free transitive action of $H^1(M, \del M;\ZZ/2)$ on the isomorphism classes of spin structures on $M$ rel. some chosen spin structure on $\del M$. 

With this in mind, we now generalize \cite[Definition 2.5]{OrsonPowell} to the case of homotopy equivalences to obtain our spin invariant. 

\begin{defn}
    Let $M$ be a $\CAT$ spin manifold with a spin structure $\fs$ and $f \maps M \to M$ an orientation preserving homotopy equivalence restricting to the identity on the boundary. Define $\Theta(f, \fs)$ to be the group element of $H^1(M, \del M;\ZZ/2)$ such that $\Theta(f, \fs)\cdot \fs=f^*(\fs)$. 
\end{defn}

Observe that \cite[Lemmas 2.6, 2.7]{OrsonPowell} apply for homotopy equivalences without change. This implies that for $M$ simply connected, $\Theta (f,\fs)$ is independent of the choice of spin structure (because $M$ has a unique spin structure up to isomorphism and isomorphic spin structures give the same $\Theta$) and also that $\Theta$ is independent of choice of using topological or smooth spin structures. For the remainder of the paper, we reflect this by denoting $\Theta(f,\fs)=:\Theta(f)$ when $M$ is simply connected. 

As in the previous section, we conclude by spelling out how we obtain an element of $H^1$ from an $h$-cobordism. 

\begin{lemma}\label[lemma]{hcob_spin_hom}
    Let $M_0, M_1$ be homeomorphic simply-connected spin $\CAT$ $4$-manifolds. Fix a homeomorphism $F\maps M_1 \to M_0$ with $F\vert_{\del M_1}=h^{-1}$ for some $h$. Then, the map $\cH^{\CAT}(M_0, M_1, h)\to H^1(M_0, \del M_0;\ZZ/2)$ given by $W \mapsto \Theta(F\circ \varphi_{W})$ is a well-defined map. Furthermore, in the case where $M_0=M_1$ and $h=\id$, we choose $F=\id$ and this map is a group homomorphism. 
\end{lemma}
\begin{proof}
    The proof is identical to that of \cref{hcob_var_hom} aside from the following changes:
    \begin{enumerate}[i)]
        \item Replacing the map $F\mapsto\Delta_F$ with $F \mapsto \Theta(F)$;
        \item Replacing \cref{htpy} with \cite[Corollary 2.8]{OrsonPowell};
        \item Replacing \cref{var_hom} with \cite[Lemma 2.7]{OrsonPowell}.
    \end{enumerate} 
    Note that the proofs from \cite{OrsonPowell} are for homeomorphisms but the same arguments work with homeomorphism replaced by homotopy equivalence and pseudoisotopy replaced by the trace of the homotopy. 
\end{proof}

\section{Surgery theory}\label[section]{surgery}
\subsection{Classical Surgery}\label[section]{classical}
In \cref{existence}, we will seek to build an $h$-cobordism realizing the invariants of \cref{thm}. To do so, we review how one builds $h$-cobordisms via classical surgery. Broadly speaking, the surgery program developed by Browder \cite{Brow72}, \cite{Brow93}, Novikov \cite{Nov}, Sullivan \cite{Sull}, and Wall \cite{Wall}, proceeds as follows: Given CAT $n$-manifolds $N, N', M$ for $n\geq 5$ with fundamental group $\pi$ (or $\CAT=\Top$, $n=4$, and $\pi$ a good group) and first Stiefel-Whitney class $w\maps \pi \to \pm 1$, 
\begin{enumerate}
    \item Find homotopy equivalences $\varphi, \varphi' \maps N, N' \to M$;
    \item Build a degree one normal bordism from $(N,\varphi)$ to $(N',\varphi')$ over $M$;
    \item Show that the associated surgery obstruction $\sigma \in L_{n+1}^h(\ZZ\pi)$ vanishes. 
\end{enumerate}
Wall packaged this information into the following exact sequence of pointed sets (cf. \cite[Chapter 10]{Wall}). 
\[\cdots L_{n+1}^h(\ZZ\pi,w_1) \dashrightarrow \cS_{\CAT}^h(M, \del M) \ra \cN_{\CAT}(M, \del M) \ra L_4^h(\ZZ\pi, w_1),\]
where the map $L^h_{n+1}(\ZZ\pi,w_1) \dashrightarrow \cS_{\CAT}^h(M, \del M)$ is a group action. We define the relevant groups as follows. The elements of the structure set $\cS_{\CAT}^h (M, \del M)$ are represented by homotopy equivalences $(\varphi,\varphi_\del) \maps (N, \del N) \to (M,\del M)$, with $N$ a $\CAT$ manifold and $\varphi_\del$ a $\CAT$-isomorphism. We mod out by the equivalence relation given by $\CAT$ $h$-cobordism over $M$, meaning that $(\varphi,\varphi_\del)\sim(\varphi',\varphi_\del')$ if there exists a $\CAT$ $h$-cobordism $(W, N, f_0, N', f_1, g)$ and a map $\Phi\maps W \to M \times [0,1]$ extending $\varphi,\varphi'$ such that $\Phi\vert_{\del_V W}\maps \del_V W \to \del M \times [0,1]$ is a $\CAT$-isomorphism and the following diagram commutes. 
\begin{equation}\label{equiv}
    \begin{tikzcd}
N' \arrow[r, "\varphi'"] \arrow[d, "f_1"]  & M \arrow[d, "\id"]        \\
\del_+W  \arrow[d, "i_+"', hook] & M\times \{1\} \arrow[d, "j_1", hook'] \\
W \arrow[r, "\Phi"]                         & M \times [0,1]                \\
\del_-W \arrow[u, "i_-", hook']  & M\times \{0\} \arrow[u, "j_0"', hook] \\
N \arrow[r, "\varphi"] \arrow[u, "f_0"]   & M \arrow[u, "\id"']    .  
\end{tikzcd}
\end{equation}
Similarly, we can define the normal invariants $\cN_{\CAT}(M,\del M)$ to be the set of degree one normal maps $N \to M$ modulo $\CAT$ bordisms equipped with a degree one normal map extending these maps. Lastly, $L_{n}^h(\ZZ\pi,w), L_{n+1}^h(\ZZ\pi,w)$ are the $L$-groups; we refer to \cite[Chapters 11 and 12]{Ranicki} for a detailed description but recall that $L_4^h(\ZZ,1)\cong \ZZ$ and $L_5(\ZZ,1)\cong 0$. When $w=1$, we often drop $w$ from the notation and simply write $L_n^h(\ZZ\pi)$. 

\begin{rmk}
    There exist a few different ways to define the structure set and normal invariants for manifolds with boundary, as studied in \cite{structuresets}. For our purposes, we use the version that they refer to as the $\reld$ case, which uses $s_\del$ and $t_\del$-triangulations (see e.g. \cite[Page 558]{ranickiexactseq} and \cite[Section 2]{structuresets}). 
\end{rmk}
The case $\pi_1=1$ is special in that the surgery exact sequence applies in the smooth category as well, via the following theorem proved by Wall (\cite[Theorem 16.5]{Wall}]) (see also corrections in \cite{CochranHabegger}). 
\begin{thm}[{\cite{Wall,CochranHabegger}}]\label[thm]{ses}
Let $M$ be a simply-connected smooth $4$-manifold. Then, the following sequence is exact.
\[0\ra \cS_{\CAT}(M, \del M) \xrightarrow{\eta} \cN_{\CAT}(M,\del M)\xrightarrow{\theta} L_4(\ZZ).\]
\end{thm}
The proof of this theorem consisted of finding explicit inverses for elements of $\ker \theta$. We will use these later in \cref{pinchy} and hence, will define them here. 

First, recall that the topological normal invariants can be identified with $[M/\del M, \G/\Top]$, and the fourth stage of a Postnikov tower for $\G/\Top$ is given by $K(\ZZ/2,2)\times K(\ZZ,4)$. Hence, 
\[[M/\del M, \G/\Top]\cong H^2(M, \del M;\ZZ/2) \oplus H^4(M, \del M).\]
Under this identification, the image of the injection $\cN_{\Diff}(M, \del M) \to \cN_{\Top}(M, \del M)$ is $\{(a,b)\vert a^2=b\}$ (see e.g. \cite[Lemma 2.1]{hollingsworth_morgan}). The surgery obstruction is given by $\theta([f\maps N \to M])=(\sigma(N)-\sigma(M))/8$ under the identification $L_4(\ZZ)\cong \ZZ$; hence, under this identification, we have 
\[\ker \theta=\{x\in H^2(M, \del M;\ZZ/2) \vert x^2=0\},\] 
which is isomorphic to $\ker w_2\subset H_2(M;\ZZ/2)$ via Poincar\'e duality when we view $w_2$ as a homomorphism $H_2(M,\ZZ/2)\to \ZZ/2$. Noting that the Hurewicz map is an isomorphism every element of $H_2(M;\ZZ/2)$ can be represented by a sphere so we define the following. 
\begin{defn}
    Given $x \in \pi_2(M)$, the \textit{Novikov pinch map} associated to $x$ is the composition 
    \[n_x \maps M \rightarrow  M \vee S^4 \xrightarrow{\id\vee (h \circ \Sigma h)} M \vee S^2 \xrightarrow{\id \vee x} M,\]
    where the first map is given by contracting the boundary of a local ball, $h$ is the Hopf map, and $\Sigma h$ is the suspension of the Hopf map. 
\end{defn}
In \cite[Theorem 5.1]{CochranHabegger} ,Cochran-Habegger computed the normal invariant of a Novikov pinch map to be \[\eta(n_x)=(1+w_2(x)) \PD^{-1}(x)\in H^2(M;\ZZ/2).\] 
In particular, for $x\in \ker w_2$, we get $\eta(n_x)=\PD^{-1}(x)$ so these maps provide inverses for $\ker \theta$. 

\subsection{Modified Surgery}\label[section]{modified}
A common struggle in applying classical surgery theory lies in finding a homotopy equivalence in the first step. To avoid this issue, Kreck developed modified surgery theory \cite{SAD}, which we will be using to prove injectivity of \cref{thm}. Morally speaking, the idea here is to restrict the homotopy and normal bundle data to the $k$-skeleton of the manifold. This translates roughly to looking at a Postnikov-Moore tower for the stable normal bundle. As in the classical setting, we can describe this in the following three-step process: for $n$-manifolds $N, N'$ as before and $k\geq \lfloor n/2 \rfloor-1$, we need to
\begin{enumerate}
    \item Find the normal $k$-type $B$ of $N,N'$ and choose $k$-smoothings $\bar \nu \maps N\to B$, $\bar \nu' \maps N'\to B$;
    \item Find a cobordism from $(N,\bar \nu)$ to $(N', \bar \nu')$ over the normal $k$-type;
    \item Show that the associated surgery obstruction $\theta \in \ell_6(\ZZ\pi,w)$ is elementary. 
\end{enumerate}

Before giving a more precise statement, we briefly recall some of the definitions involved; for further details, we refer to \cite{SAD} for the overall theory and to \cite{SixtMonoids} for a detailed study of the $\ell_6$ monoid. Also, we remark that as in the classical setting, we can work with based modules to replace $h$-cobordism with $s$-cobordism but we omit this as it is unnecessary for our purposes. 

To start, let $B$ be a space with the homotopy type of a CW-complex and $\xi \maps B \to \BCAT$ a fibration. An \textit{$n$-dimensional $(B, \xi)$-$\CAT$-manifold} is a pair $(M ,\bar \nu)$ where $M$ is an $n$-dimensional $\CAT$-manifold and $\bar \nu$ is a lift of the stable normal bundle. When the category doesn't matter, we omit $\CAT$ and say $(B,\xi)$-manifold. A \textit{$(B, \xi)$-nullbordism} of a $(B,\xi)$-manifold $(M, \bar \nu)$ is a $(n+1)$-dimensional $(B,\xi)$-manifold $(W, \bar\nu_W)$ such that $(\del W, \bar \nu_W\vert_{\del W})=(M, \bar \nu)$. A \textit{$(B,\xi)$-cobordism} between two $n$-dimensional $(B, \xi)$-manifolds $(M, \bar \nu_{M}), (N, \bar \nu_N)$ with a fixed boundary identification $f\maps \del M \to \del N$ is a $(B,\xi)$-nullbordism of $(X, \bar\nu_X)$, where $X:=M \cup_f -N$  and $\bar\nu_X$ is a lift of the stable normal bundle of $X$ such that 
\begin{enumerate}[i)]
    \item $\bar\nu_X\vert_{M} \simeq \bar \nu_M$ and $\bar\nu_X\vert_{N} \simeq \bar \nu_N$;
    \item $\bar\nu_X\vert_{\del M}=\bar\nu_X\vert_{\del N}\circ f \maps \del M \to B$.
\end{enumerate}
With this in hand, we now define the normal $k$-type. 

\begin{defn}
Let $B$ be a space with the homotopy type of a CW-complex with a finite $(k+2)$-skeleton and $\xi \maps B \to \BSCAT$ a fibration. For $(M, \bar \nu)$ a $(B, \xi)$-manifold, we say that $(M,\bar \nu)$ is a \textit{normal $k$-smoothing} into $(B,\xi)$ if $\bar \nu$ is $(k+1)$-connected. If $\xi$ is $(k+1)$-coconnected and there exists a normal $k$-smoothing $(M, \bar \nu)$ into $(B,\xi)$, then we call $(B, \xi)$ the \textit{normal $k$-type} of $M$. 
\end{defn}

Now, we proceed to briefly review the surgery obstruction and the even dimensional $\ell$ monoid. Let $\pi$ be a group with a homomorphism $w \maps \pi \to \ZZ/2$. The involution on $\ZZ\pi$ defined by sending $g\in \pi$ to $\bar g= w(g)g^{-1}$ allows us to consider right $\ZZ\pi$-modules as left $\ZZ\pi$-modules via the action $g\cdot x= x\bar g$ (in particular giving $V^*$ the structure of a left $\ZZ\pi$-module). Recall that a \textit{$(-1)^q$-quadratic form} $(V, \lambda, \mu)$ over $\ZZ\pi$ consists of a left $\ZZ\pi$-module $V$ equipped with a $(-1)^q$-hermitian form $\lambda$ and a quadratic refinement $\mu \maps V \to \ZZ\pi/\<x-(-1)^q\bar x\>$. We can also define a reduced $(-1)^q$-quadratic form $(V, \lambda, \tilde \mu)$ by replacing $\mu$ with $\tilde \mu\maps V \to \ZZ\pi/\<x-(-1)^q\bar x, 1\>$. With these in hand, we can define the objects that make up the even dimensional $\ell$ monoid. 

\begin{defn}[{\cite[Section 4]{SAD},\cite[Definition 2.4]{SixtMonoids}}]\label[defn]{preformation}
A \textit{$(-1)^q$-preformation} is a triple
\[(V_0 \xleftarrow{f} V \xrightarrow{g} V_1, \lambda, \mu),\]
such that 
\begin{enumerate}[i)]
    \item V is a finitely generated free left $\ZZ\pi$-module and $V_0, V_1$ are free.
    \item $\lambda\maps V_0 \to V_1^*$ is an isomorphism.
    \item $\lambda_V (x,y):= \lambda(f(x), g(y))$ and $\mu$ comprise a $(-1)^q$-quadratic form on $V$. 
\end{enumerate}
\end{defn}
\begin{rmk}
    In \cite{SixtMonoids}, Sixt refers to $(-1)^q$-preformations as regular split $(-1)^q$-preformations to remain consistent with definitions in algebraic surgery. 
\end{rmk}
We can think of a $(-1)^q$-quadratic form $(V,\lambda,\mu)$ as a $(-1)^q$-preformation by setting $V_0=V_1=V$ and $f=g=\id$. In particular, consider the $(-1)^q$-hyperbolic form, which is the quadratic form $(\ZZ\pi\oplus \ZZ\pi, H, \mu)$, where
\[H=\begin{pmatrix}
    0& 1\\
    (-1)^q&0
\end{pmatrix},\]
and $\mu$ is given by multiplication, i.e. $\mu(x,y)=xy$. We want to look at preformations up to stabilization by hyperbolics, which leads to the following definition.
\begin{defn}[{\cite[Section 4]{SAD}}]
    The monoid $\ell_{2q}(\ZZ\pi,w)$ is the set of $(-1)^q$-preformations under the direct sum operation modulo the relation $z \sim z'$ if $z \oplus H^{\oplus n} \cong z' \oplus H^{\oplus m}$ for some $n,m\in \NN$. Working with reduced quadratic forms instead, the same definition gives $\tilde \ell_{2q}(\ZZ\pi,w)$. 
\end{defn}

We will now recall how to obtain the surgery obstruction associated to a $(B,\xi)$-bordism $(Y,\bar \nu_Y)$ between $(2q-1)$-dimensional normal $(q-2)$-smoothings $(M,\bar\nu_M)$ and $(N,\bar\nu_N)$. Note that we require $q\geq 3$ because surgery only works in dimension at least $4$. This was done in \cite[Section 5]{SAD} but we repeat the main points here for completeness. We will start off by suggesting (without justification) that
\[H_q(Y, M;\ZZ\pi)\xleftarrow{f}\im(d\maps \pi_{q+1}(B,Y)\to \pi_q(Y))\xrightarrow{g}H_q(Y, N;\ZZ\pi),\]
is a natural starting point for building the surgery obstruction, where $f$ is the composition of the inclusion $H_q(Y;\ZZ\pi)\to H_q(Y,M;\ZZ\pi)$ with the Hurewicz map $\pi_q(Y)\to H_q(Y;\ZZ\pi)$ and similarly for $g$. One immediate worry is that condition i) of \cref{preformation} might not hold. So, we first have to modify $Y$ to ensure that it does hold. By \cite[Proposition 4]{SAD}, we can assume $Y\to B$ is $q$-connected via surgery below the middle dimension. Then, denote $Y_n$ to be $Y\#n(S^q \times S^q)$. Note that there exists a standard $B$-structure on $S^q \times S^q$ by restricting the unique $B$-structure on $D^{q+1}$ to $S^q$. 
\begin{lemma}\label[lemma]{free}
    For some $n\geq 0$, the $\ZZ\pi$-modules $H_q(Y_n,M;\ZZ\pi)$ and $H_q(Y_n,N;\ZZ\pi)$ are finitely generated free $\ZZ\pi$-modules. 
\end{lemma}
\begin{proof}
    Then, for $i<q$, the map $\pi_i(M)\twoheadrightarrow \pi_i(B)$ factors through $\pi_i(Y)\cong \pi_i(B)$ via the inclusion so $\pi_i(M)\twoheadrightarrow \pi_i(Y)$. Thus, we  obtain $\pi_i(Y,M)\cong 0$ and hence, $H_i(Y,M;\ZZ\pi)=0$, with the same argument also holding for $N$. The universal coefficient spectral sequence (UCSS) gives us that for any $\ZZ\pi$-module $\Lambda$, $H_i(Y,N;\Lambda)=0$. For $i>q$, the same is true by Poincar\'e duality and the UCSS. Hence, by \cite[Lemma 2.3]{Wall}, $H_q(Y,M;\ZZ\pi)$ and $H_q(Y,N;\ZZ\pi)$ are finitely generated stably free modules, so for some $n$, $H_q(Y_n,M;\ZZ\pi)$ and $H_q(Y_n,N;\ZZ\pi)$ are finitely generated free $\ZZ\pi$-modules. 
\end{proof}
We also note that $\im(d\maps \pi_{q+1}(B,Y_n)\to \pi_q(Y_n))$ is finitely generated, e.g. via \cite[Lemma 3.4]{CCPS}. 

We now proceed to describe $\lambda$ and $\mu$. In the proof of \cref{free}, we showed that for $i<q$, $H_i(Y,N;\ZZ\pi)=0$. The same proof applies to $H_i(Y_n,N;\ZZ\pi)$ as well; hence, the evaluation map on $ H^q(Y_n,N;\ZZ\pi)$ is an isomorphism via the UCSS. Hence, $\lambda:=\ev\circ\PD^{-1}\maps H_q(Y_n,M;\ZZ\pi)\to H_q(Y_n,N;\ZZ\pi)^*$ is an isomorphism.

We now describe the required quadratic form. Denote the group of regular homotopy classes of immersed $q$-spheres in $Y$ via $\imm_q(Y)$. In \cite[Theorem 5.2]{Wall}, Wall defined a pair of functions that we will denote as $\lambda^{\geom}\maps\imm_q(Y_n)\times\imm_q(Y_n)\to \ZZ\pi$ and $\mu\maps \imm_q(Y_n)\to \ZZ\pi/\<x-(-1)^qx\>$ (see also \cite[Definition 11.21]{Ranicki}). Roughly speaking, these count the intersections and self-intersections of spheres. When restricted to immersed spheres with stably trivial normal bundle, this becomes a quadratic form. After possibly summing on a $S^q\times S^q$ to $Y_n$, Kreck showed in \cite[Proposition 6]{SAD} that for $q\neq 1,3,7$, every element of $\im(d\maps \pi_{q+1}(B,Y_n)\to \pi_q(Y_n))$ can be represented by a unique such sphere up to homotopy. Hence, we obtain a quadratic form on $\im(d\maps \pi_{q+1}(B,Y_n)\to \pi_q(Y_n))$, which by an abuse of notation, we also write as $\lambda^{\geom}$. In the case where $q=3,7$, elements of $\im(d\maps \pi_{q+1}(B,Y)\to \pi_q(Y_n))$ fail to have unique representatives. In this scenario, we can still define $\lambda^{\geom}$ in the same way because it is a homotopy invariant (follows from e.g. \cite[Proposition 7.22]{Ranicki}) but $\mu$ may not be well defined. When $\<\pi_{q+1}(B), w_{q+1}(B)\> \neq 0$, Kreck showed that any two representatives differ under $\mu$ by a multiple of $1$. Hence, replacing $\mu$ with the reduced self-intersection number $\tilde \mu \maps \imm_q(Y_n)\to \ZZ\pi/\<x-(-1)^qx,1\>$ provides a quadratic form. To address the case where $q=3,7$ and $\<\pi_{q+1}(B), w_{q+1}(B)\> = 0$, Kreck constructed a $\mu$ on $\im(d\maps \pi_{q+1}(B,Y_n)\to \pi_q(Y_n))$ giving a quadratic refinement of $\mu$ in \cite[Proposition 7]{SAD}. 

It remains to check that the map $\lambda$ induces on $\im(d\maps \pi_{q+1}(B,Y_n)\to \pi_q(Y_n))$ is equal to $\lambda^{\geom}$ up to sign. Define $f_n\maps \im(d\maps \pi_{q+1}(B,Y_n)\to \pi_q(Y_n))\to H_q(Y_n,M;\ZZ\pi)$ as before and likewise for $g_n$.  

\begin{lemma}
    Let $x,y\in\im(d\maps \pi_{q+1}(B,Y_n)\to \pi_q(Y_n))$. Then, $\lambda(f_{n}(x),g_{n}(y))=\lambda^{\geom}(x,y)$. 
\end{lemma}
\begin{proof}
    It is well known that the form $\lambda^{\geom}$ induces on homology is the same as the algebraic intersection form $\lambda^{\alg}$, e.g. via \cite[Proposition 7.22]{Ranicki}. In other words, the following diagram commutes:

    \[\begin{tikzcd}
\imm_q(Y_n)\times \imm_q(Y_n) \arrow[r, "\lambda^{\geom}"] \arrow[d]      & \ZZ\pi \\
\pi_q(Y_n)\times \pi_q(Y_n)                                     &        \\
\arrow[u,"\cong"] \pi_q(\tilde Y_n)\times \pi_q(\tilde Y_n) \arrow[d]                       &        \\
H_q(Y_n;\ZZ\pi)\times H_q(Y_n;\ZZ\pi) \arrow[ruuu, "\lambda^{\alg}", bend right] &       
\end{tikzcd}\]

We claim that the following diagram also commutes:
\[\begin{tikzcd}
H_q(Y_n;\ZZ\pi)\times H_q(Y_n;\ZZ\pi) \arrow[r, "\lambda^{\alg}"] \arrow[d, "(f_n)_* \times (g_n)_*"'] & \ZZ\pi \\
{H_q(Y_n,M;\ZZ\pi)\times H_q(Y_n,N;\ZZ\pi)} \arrow[ru, "\lambda"]                              &       
\end{tikzcd}\]

In particular, if this diagram commutes, then the combination with the first diagram will prove the statement of the lemma, as the vertical composition starting at $\pi_q(Y_n)\times \pi_q(Y_n)$ sends $(x,y)\mapsto (f_n(x),g_n(y))$

We prove commutativity of the diagram by proving that the adjoint of $\lambda^{\alg}$ agrees with $g^*\circ \lambda\circ f_*$. Note that the evaluation map is natural so it commutes with $g^*$; hence, $g^*\circ \lambda \circ f=\ev \circ g^* \circ \PD^{-1} \circ f$. So, the following diagram commutes, where the bottom row is the adjoint of $\lambda^{\alg}$, proving the desired equality.
\[\begin{tikzcd}
H_q(Y_n;\ZZ\pi) \arrow[d, "\id"] \arrow[r, "f"] & {H_q(Y_n,M;\ZZ\pi)} \arrow[d] \arrow[r, "\PD^{-1}"] & {H^q(Y_n,N;\ZZ\pi)} \arrow[r, "g^*"] \arrow[d, "g^*"'] & H^q(Y_n;\ZZ\pi) \arrow[r, "\ev"] \arrow[d, "\id"'] & H_q(Y_n;\ZZ\pi)^* \arrow[d, "\id"] \\
H_q(Y_n;\ZZ\pi) \arrow[r, "i_*"]                & {H_q(Y_n,M \cup N;\ZZ\pi)} \arrow[r, "\PD^{-1}"]   & H^q(Y_n;\ZZ\pi) \arrow[r, "\id"]                        & H^q(Y_n;\ZZ\pi) \arrow[r, "\ev"]                   & H_q(Y_n;\ZZ\pi)^*                 
.\end{tikzcd}\]
\end{proof}

In particular, this implies that $\lambda$ and $\mu$ provide a $(-1)^q$-quadratic form on $V$. With this in mind, we can define the surgery obstruction. 
\begin{defn}[{\cite[Section 4]{SAD}}]
    Let $(Y,\bar\nu_Y)$ be a $(B,\xi)$-bordism between two $(2q-1)$-dimensional normal $(q-2)$-smoothings $(M,\bar\nu_M)$ and $(N,\bar\nu_N)$ for $q\geq 3$. Surger $Y$ so that $Y\ra B$ is $q$-connected and connect sum some number of $S^q \times S^q$s to obtain $Y_n$ such that $H_q(Y_n,M;\ZZ\pi)$, $H_q(Y_n,N;\ZZ\pi)$ are free $\ZZ\pi$-modules.
    
    When $q\neq 3,7$ or $q=3,7$ and $\<\pi_{q+1}(B),w_{q+1}(B)\>= 0$, the \textit{surgery obstruction} associated to $(Y,\bar\nu_Y)$ is 
    \[\theta(Y,\bar\nu_Y):=[H_q(Y_n, M;\ZZ\pi)\xleftarrow{f}\im(d\maps \pi_{q+1}(B,Y_n)\to \pi_q(Y_n))\xrightarrow{g}H_q(Y_n, N;\ZZ\pi), \lambda, \mu],\]
    where $\lambda\maps H_q(Y_n,M;\ZZ\pi)\to H_q(Y_n,N;\ZZ\pi)^*$ and $\mu\maps \im(d\maps \pi_{q+1}(B,Y_n)\to \pi_q(Y_n)) \to \ZZ\pi/\<x-(-1)^qx\>$ are as above, and the maps $f,g$ are induced by the inclusion composed with the Hurewicz map. 
    
    When $q=3,7$ and $\<\pi_{q+1}(B),w_{q+1}(B)\>\neq 0$, the same definition with $\mu$ replaced by the reduced Wall self-intersection form $\tilde \mu\maps \im(d\maps \pi_{q+1}(B,Y_n)\to \pi_q(Y_n)) \to \ZZ\pi/\<x-(-1)^qx,1\>$ gives a surgery obstruction in $\tilde \ell_{2q}(\ZZ\pi,w)$ instead. 
\end{defn}
We note that the surgery obstruction is well defined because any two choices of $Y_n$ are stably diffeomorphic rel. boundary over $(B, \xi)$ via \cite[Corollary 3]{SAD}, and hence, determine the same elements of the monoid. 

As outlined above, we care about deciding when the surgery obstruction is elementary:
\begin{defn}[{\cite[Section 4]{SAD}}]\label[defn]{elementary}
    A $(-1)$-preformation $(V_0 \xleftarrow{f} V \xrightarrow{g} V_1, \lambda, \mu)$ is \textit{elementary} if there exists a free submodule $U\subset V$ such that
\begin{enumerate}[i)]
    \item $f,g$ map $U$ injectively to direct summands $U_i\subset V_i$.
    \item $\proj_{U_1^*}\circ \lambda\vert_{U_0}\maps U_0 \to U_1^*$ is the zero map and $ \mu\vert_U=0$.
    \item $\proj_{(V_1/U_1)^*}\circ \lambda\vert_{U_0} \maps U_0\to (V_1/U_1)^*$ is an isomorphism. 
\end{enumerate}
Given an element $x\in \ell_{2q}(\ZZ\pi_1(B),w_1(B))$, we say that it is elementary if it admits an elementary representative. The same definitions hold for $\tilde \ell_{2q}(\ZZ\pi_1(B),w_1(B))$ with $\mu$ replaced by $\tilde \mu$. 
\end{defn}
In both ii) and iii) we note that $V_1^*$ splits as $U_1^* \oplus (V_1/U_1)^*$ by i) so we compose with the projection is well defined. Now that we have defined the relevant terms, we can state Kreck's main theorem for even-dimensional surgery

\begin{thm}[{\cite[Theorem 3]{SAD}}]\label[thm]{mod_surg}
Let $(M,\bar\nu_M)$ and $(N,\bar\nu_N)$ be $(2q-1)$-dimensional normal $(q-2)$-smoothings, $q\geq 3$. Additionally suppose that $q \neq 3,7$ or $\<w_{q+1}(B), \pi_{q+1}(B)\>=0$. Let $(Y,\bar\nu_Y)$ be a $(B,\xi)$-bordism between $M$ and $N$. Then, $(Y,\bar\nu_Y)$ is $(B,\xi)$-bordant rel. boundary to a relative $s$-cobordism if and only if the associated surgery obstruction $\theta(Y, \bar \nu)\in \ell_{2q}(\ZZ\pi_1(B),w_1(B))$ is elementary. If $q=3,7$ and $\<w_{n+1}(B), \pi_{n+1}(B)\>\neq 0$, then the same statement holds with the reduced surgery obstruction $\theta(Y, \bar\nu) \in \tilde \ell_{2n}(\ZZ\pi_1(B),w_1(B))$. 
\end{thm}

We will be using this theorem to construct an $h$-cobordism between simply connected $5$-manifolds; hence, we now specialize to the case $\pi=1$ and $q=3$. For the rest of this section, assume that all $(B,\xi)$ nullbordisms ($Y, \bar \nu_Y)$ have been surgered so that $Y\to B$ is $3$-connected. Note that in this scenario, $H_n(Y,N)\cong H_n(Y,M)^*$ so $H_n(Y,N)$ is already free without needing to stabilize and likewise for $H_n(Y,M)$. Hence, we can use $Y_n=Y$ in the definition of the surgery obstruction. 

For a simply-connected CAT manifold $M$, the normal $1$-type $B\xrightarrow{\xi}\BSCAT$ is given by $\BSCAT\xrightarrow{\id}\BSCAT$ if $M$ is not spin and $\BCATSpin\xrightarrow{p}\BSCAT$ if $M$ is spin, where $p$ is induced by the covering map (see e.g. \cite[Lemmas 3.1, 3.5]{KPLT}). We start by determining that $\<w_{4}(B), \pi_{4}(B)\>$ is nontrivial. 

\begin{lemma}\label[lemma]{w4}
    Let $B\in\{\BSCAT,\BCATSpin\}$. Then, $\<w_{4}(B), \pi_{4}(B)\>\neq 0$. 
\end{lemma}
\begin{proof}
    Before we start, we unpack some notation. First of all, $w_{4}(\BSCAT)$ refers to the $w_4(\gamma)$, where $\gamma$ is the universal bundle over $\BSCAT$. Likewise, $w_4(\BCATSpin)=w_4(p^*\gamma)$. Also, when we write the evaluation $\<w_{4}(B), \pi_{4}(B)\>$, we view elements of $\pi_4(B)$ as homology elements via the Hurewicz map. 

    Now, we start by showing $\<w_{4}(\BSCAT), \pi_{4}(\BSCAT)\>\neq 0$. Note that (e.g. via \cite[Page 93]{Hatchervb}) there exists a rank $4$ vector bundle $E$ over $S^4$ such that the Euler class $e(E)$ generates $H^4(S^4)$, namely the bundle associated to the quaternionic Hopf algebra. Let $f \maps S^4 \to \BSCAT$ be the classifying map for this bundle. Then, by the naturality of the evaluation map,
    \begin{align*}
        \<w_4(\BSCAT),[f]\>&=\<f^*w_4(\gamma),[S^4]\>\\
        &=\<w_4(f^*\gamma),[S^4]\>\\
        &=\<w_4(E),[S^4]\>\\
        &=\< \red_2 e(E), [S^4]\>\\
        &=\< 1, [S^4]\>\\
        &=1.
    \end{align*}
    Hence, $\<w_{4}(\BSCAT), \pi_{4}(\BSCAT)\>\neq 0$.

    Now, considering $\BCATSpin$, we note that $p_*$ is an isomorphism on $\pi_4$ so and commutes with the Hurewicz map so
    \begin{align*}
        \<w_{4}(\BCATSpin), p_*^{-1}([f])\>=\<w_{4}(p^*\gamma), p_*^{-1}([f])\>=\<p^*w_{4}(\gamma), p_*^{-1}([f])\>=\<w_{4}(\gamma), [f]\>=1.
    \end{align*}
\end{proof}

As a result, the surgery obstruction lies in $\tilde \ell_6(\ZZ)$. In particular, $\ZZ\pi/\<x-\bar x,1\>=\ZZ/\<1\>=0$ so we can ignore quadratic refinements. Thus, we omit $\tilde \mu$ when writing preformations in $\tilde \ell_6(\ZZ)$. It will also be convenient to write our surgery obstruction entirely in terms of homology, which we can do in this case.

\begin{lemma}
    The surgery obstruction $\theta$ associated to a $(B,\xi)$-bordism $(Y,\bar\nu_Y)$ between two simply-connected $5$-dimensional normal $1$-smoothings $(M,\bar\nu_M)$ and $(N,\bar\nu_N)$ satisfies the equality
    \[\theta=[H_3(Y, M)\xleftarrow{f'}H_3(Y)\xrightarrow{g'}H_3(Y, N), \lambda]\in \tilde \ell_6(\ZZ),\]
    where $\lambda$ is as before and $f',g'$ are induced by the inclusion. 
\end{lemma}
\begin{proof}
    Recall that $f$ and $g$ both factor through the Hurewicz map so we obtain the following commutative diagram, where $f',g'$ are induced by the inclusion: 
\[\begin{tikzcd}
{H_3(Y,M)} & {\im(d\maps \pi_4(B,Y)\to \pi_3(Y))} \arrow[r, "g"] \arrow[l, "f"'] \arrow[d, "\cong"] & {H_3(Y,N).} \\
           & H_3(Y) \arrow[ru, "g'"'] \arrow[lu, "f'"]                                                &            
\end{tikzcd}\]
To see that the vertical map is an isomorphism, first observe that $\pi_3 (B)=0$ so $\pi_4(B,Y)$ surjects onto $\pi_3(Y)$ and thus, $\im(d\maps \pi_4(B,Y)\to \pi_3(Y))=\pi_3(Y)$. Then, the Hurewicz map $h\maps \pi_3(Y)\to H_3(Y)$ is a surjection because $Y$ is simply connected. For injectivity, observe that the zeroth, first, and second rational homotopy groups of $Y$ vanish so by the rational Hurewicz theorem, $h\otimes \QQ \maps \pi_3(Y)\otimes \QQ \to H_3(Y,\QQ)$ is an isomorphism. In particular, localization preserves kernels so $\ker h \otimes \QQ=0$ implying that $\ker h$ must be torsion, but $\pi_3(Y)$ was torsion-free, so $\ker h=0$.  

As a result, the diagram provides an isomorphism of preformations, so the surgery obstruction is equal to the desired preformation. 
\end{proof}
\begin{rmk}
    By \cite[Lemma 2.9]{SixtMonoids}, it would have been enough to show that the vertical map is a surjection. 
\end{rmk}
For the sake of notation, we redefine $f$ and $g$ to be the $f'$ and $g'$ from the above lemma. 

We can now simplify the process of checking elementariness to the following simpler condition, as in \cite{Kreck_hcob}. 
\begin{lemma}\label[lemma]{sc_elem}
    The surgery obstruction $\theta$ associated to a $(B,\xi)$-bordism $(Y,\bar\nu_Y)$ between two simply-connected $5$-dimensional normal $1$-smoothings $(M,\bar\nu_M)$ and $(N,\bar\nu_N)$ is elementary if there exists a free subgroup $U\subset H_3(Y)$ such that
\begin{enumerate}[i)]
    \item The inclusion maps $U$ injectively to direct summands $U_0, U_1$ in $H_3(Y,M)$ and $H_3(Y,N)$;
    \item The intersection form vanishes on $U$;
    \item $\rk U$=$\rk H_3(Y, M)/2$.
\end{enumerate}
\end{lemma}
\begin{proof}
We will compare conditions i)-iii) of the lemma with those of \cref{elementary}. Condition i) is unchanged so we start with condition ii). As stated above, the bilinear form $\lambda$ induced on $H_3(Y)$ is the intersection form. Thus, it is enough to show that the intersection form vanishes on $U$. 

For condition iii), condition ii) gives us that under the splitting $H_3(Y,N)^*\cong U_1^*\oplus (H_3(Y,N)/U_1)^*$, the image of $U_0$ under $\lambda$ lies in the $(H_3(Y,N)/U_1)^*$ summand. This is an equality if and only if we have $\rk \lambda(U_0)=\rk (H_3(Y,N)/U_1)^*$ because $\lambda$ was an isomorphism and $U_0$ is a direct summand of $H_3(Y,M)$. Noting that $\rk (H_3(Y,N)/U_1)^*=\rk H_3(Y,N) - \rk U_1$ and $\rk U = \rk U_0=\rk U_1$, we arrive at the desired rank condition. 
\end{proof}

Observe that $\ker f \subset \ker (H_3(Y)\to H_3(Y, \del Y))$ because the map $H_3(Y)\to H_3(Y,\del Y)$ factors through $H_3(Y,M)$. Hence, $\lambda$ descends to a well-defined map on $H_3(Y)/\ker f$. As a result, when further conditions are satisfied, we can ignore the direct sum condition by working in $H_3(Y)/\ker f$ instead. 

\begin{lemma}\label[lemma]{free_elem}
Let $(Y,\bar\nu_Y)$ be a $(B,\xi)$-bordism between simply-connected $5$-dimensional normal $1$-smoothings $(M,\bar\nu_M)$ and $(N,\bar\nu_N)$ such that $H_2(M)$ and $H_2(N)$ are free and $\ker f = \ker g$. Then, $\theta$ is elementary if there exists $U\subset H_3(Y)/\ker f$ such that the intersection form vanishes on $U$ and $\rk U$=$\rk H_3(Y, M)/2$.
\end{lemma}
\begin{proof}
    Given such a $U$, we will construct $U'$ as in \cref{sc_elem}. 
    
    Consider the exact sequence
    \[0 \ra H_3(M)\xrightarrow{i} H_3(Y)\xrightarrow{f} H_3(Y,M),\]
    where we recall that $H_4(Y,M)=0$ from the proof of \cref{free}. Note that $i$ is a split injection because $H_3(Y,M)$ is free so let $s$ be a splitting $H_3(Y)/\ker f \ra H_3(Y)$. Consider the subgroup $V\subset H_3(Y)/\ker f$ defined via 
    \[V=\{x \in H_3(Y)\vert nx\in U\ \text{for some $n\in \ZZ$}\},\]
    which is a direct summand of $H_3(Y)/\ker f$ because $\ZZ$ is a PID. Define $U'$ to be $s(V)$.

    Note that the long exact sequence of a pair gives us the short exact sequence
    \[0 \ra H_3(Y)/\ker f \xrightarrow{\bar f} H_3(Y,M) \ra \im(H_3(Y,M)\ra H_2(M)) \ra 0,\]
    where $\bar f$ is given by $f$ descending to the quotient. This splits because $H_2(M)$ is free. Note that $f$ factors through $H_3(Y)/\ker f$ so $f(U')=\bar f(V)$, which is a direct summand in $H_3(Y,M)$ because $V$ was a direct summand of $H_3(Y)/\ker f$ and $\bar f$ is a split injection. By assumption, $\ker f=\ker g$ so $H_3(Y)/\ker f=H_3(Y)/\ker g$ and the same argument shows that $U'$ maps to a direct summand in $H_3(Y,N)$ as well. Hence, condition i) is satisfied. 

    By bilinearity, the intersection form vanishes on $V$ because it does on $U$. So, the intersection form must vanish over $U'$ and condition ii) holds as well. 

    Lastly, observe that $V$ has same rank as $U$ so $\rk U'=\rk V=\rk H_3(Y,M)/2$ and condition iii) is satisfied by $U'$. 
\end{proof}

We finish this section by constructing a candidate for $U$ which reduces the process of checking elementariness to a homology rank computation which can be done, for instance, using Mayer-Vietoris (as we will do in \cref{uniqueness}). 

\begin{prop}\label[prop]{rk_elem}
    Let $(Y,\bar\nu_Y)$ be a $(B,\xi)$-bordism between simply-connected $5$-dimensional normal $1$-smoothings $(M,\bar\nu_M)$ and $(N,\bar\nu_N)$ such that $H_2(M)$ and $H_2(N)$ are free and $\ker f = \ker g$. Denote $X=\del Y=M\cup_h -N$, where $h$ comes from the identification of $\del M$ with $\del N$. Then, the surgery obstruction is elementary if $\rk H_3(X)=\rk H_2(M)+\rk H_3(M)$. 
\end{prop}
\begin{proof}
    To start, recall that $H_4(Y,M)\cong 0$ so we have an injection $H_3(M)\hookrightarrow H_3(Y)$. Additionally, $H_4(Y,X)\cong H^2(Y)$, which is trivial by the UCT so $H_3(X)\hookrightarrow H_3(Y)$, which implies that we have the chain of injections $H_3(M)\hookrightarrow H_3(X)\hookrightarrow H_3(Y)$. To simplify notation, we will identify these groups with their images throughout this proof. 
    
    We need to find $U$ as in \cref{free_elem}. Note that $\lambda$ is skew symmetric and unimodular (and hence, symplectic) on $H_3(Y)/H_3(X)\subset H_3(Y)/H_3(M)$. So, there exists a half rank direct summand $U_1$ on which it vanishes. Furthermore, $\lambda$ vanishes on $H_3(X)$ so in particular, it vanishes on $U_2:=H_3(X)/H_3(M)$. Hence, we can take $U:=U_1 \oplus U_2$. 
    
    By construction, $\lambda$ vanishes on $U$ so it remains to check that $\rk(U)$=$\rk H_3(Y, M)/2$. We can compute from the definition of $U$ that
    \[2\rk U= (\rk H_3(Y)-\rk H_3(X)) +2(\rk H_3(X)-\rk H_3(M))= \rk H_3(Y)+\rk H_3(X)-2 H_3(M).\] 
    From the long exact sequence of a pair, we get the exact sequence
    \[0\ra H_3(M) \ra H_3(Y)\ra H_3(Y,M) \ra \im (H_3(Y,M)\ra H_2(M)) \ra 0.\]
    Observe that \[\im (H_3(Y,M)\ra H_2(M))=\ker (H_2(M) \ra H_2(Y)),\] 
    and $H_2(Y)$ is torsion, so \[\rk \ker (H_2(M) \ra H_2(Y))=\rk H_2(W).\] 
    Hence, \[\rk H_3(Y,M)=\rk H_3(Y)+\rk H_2(M)-\rk H_3(M).\] 
    Thus, $\rk U = \rk H_3(Y,M)/2$ precisely when the conditions of the lemma are satisfied. 
\end{proof}
\section{Realizing Invariants}\label[section]{existence}
The goal of this section is to show that every variation and element of $H^1(M_0, \del M_0;\ZZ/2)$ can be realized by an $h$-cobordism, thus proving the surjectivity of the main theorem:

\begin{thm}\label[thm]{surj}
    Let $M_0, M_1$ be homeomorphic simply-connected $\CAT$ $4$-manifolds. Fix a homeomorphism $F\maps M_1 \to M_0$ with $F\vert_{\del M_1}=h^{-1}$  for some $h$ and $\Delta \in \cV(H_2(M_0), \lambda_{M_0})$. Then, there exists a $\CAT$ $h$-cobordism $(W,M_0, f_0, M_1, f_1,g)$ rel. $h$ such that $\Delta_{F\circ \varphi_W}=\Delta$. Furthermore, if $M_0$ is spin, then we can also choose $W$ so that $\Theta(\varphi_W)=x$ for any $x\in H^1(M_0, \del M_0;\ZZ/2)$.
\end{thm}

We start by proving the result in the topological category. 

\begin{prop}\label[prop]{top_surj}
    Let $M_0, M_1$ be homeomorphic simply-connected topological $4$-manifolds. Fix a homeomorphism $F\maps M_1 \to M_0$ with $F\vert_{\del M_1}=h^{-1}$  for some $h$ and $\Delta \in \cV(H_2(M_0), \lambda_{M_0})$. Then, there exists a topological $h$-cobordism $(W,M_0, f_0, M_1, f_1,g)$ rel. $h$ such that $\Delta_{F\circ \varphi_W}=\Delta$. Furthermore, if $M_0$ is spin, then we can also choose $W$ so that $\Theta(\varphi_W)=x$ for any $x\in H^1(M_0, \del M_0;\ZZ/2)$.
\end{prop}
\begin{proof}
    Fix a variation $\Delta\in \cV(H_2(M_0),\lambda_{M_{0}})$. By \cite[Theorem 3.8]{OrsonPowell}, there exists a homeomorphism $G\maps M_0 \to M_0$ restricting to the identity on $\del M_0$ such that $\Delta_G=\Delta$. Hence, consider the cobordism $W$ given by $(M_0\times [0,1],M_0, G, M_1, F, \id)$. This is an $h$-cobordism rel. $h$. Furthermore, $\varphi_W\simeq F^{-1}\circ G$ so $\Delta_{F\circ\varphi_{W}}=\Delta_G=\Delta$ by \cref{htpy} and \cref{var_hom}. 

    In the case where $M$ is spin, let $x\in H^1(M_0, \del M_0;\ZZ/2)$. By \cite[Theorem 3.8]{OrsonPowell} there exists a choice of $G$ such that $\Theta(G)=x$. Assuming that we chose such a $G$, 
    \[\Theta(F\circ \varphi_W)=\Theta(F\circ F^{-1}\circ G)=\Theta(G)=x.\]
\end{proof}

Morally speaking, the argument in the smooth category is the same. However, we don't have a realization result for diffeomorphisms \`a la \cite[Theorem 3.8]{OrsonPowell}, and hence are unable to use a product cobordism. Luckily, we can instead use Novikov pinch maps and the surgery exact sequence to achieve the same outcome. We start by showing that we can modify any pair of homotopy equivalences to be $h$-cobordant, without changing their associated variation maps. 

\begin{lemma}\label[lemma]{pinchy}
    Let $M$ be a simply connected smooth $4$-manifold with boundary. Given $f\maps N_0 \to M$, $g\maps N_1\to M$ homotopy equivalences inducing diffeomorphisms on the boundary, there exists a homotopy equivalence $n\maps M\to M$ inducing the identity on the boundary such that $[n \circ f]=[g]\in \cS^{\Diff}(M, \del M)$. Furthermore $\Delta_{n}=\Delta_{\id}$. 
\end{lemma}
\begin{proof}
    Recall from \cref{classical} that the surgery exact sequence in this case is
    \[0\ra \cS^{\Diff}(M, \del M) \xrightarrow{\eta} \cN^{\Diff}(M, \del M) \xrightarrow{\theta} L_4(\ZZ),\]
    with inverses of $\ker \theta$ being given by Novikov pinch maps, recalling that the image of $\cN^{\Diff}(M, \del M)$ in $H^2(M, \del M;\ZZ/2)$ is $\ker w_2$. Hence, let $n$ be the pinch map associated to $\eta(g)-\eta(f)$, viewed as elements of $H^2(M,\del M;\ZZ/2)$. By definition, $n$ restricts to the identity on $\del M$. Also, by \cite[Proposition 2.2]{Brumfiel}, \[\eta(n \circ f)=\eta(n)+n^*\eta(f)=\eta(n)+\eta(f)=\eta(g),\]
    noting that $n$ induces the identity on $H^2(M, \del M;\ZZ/2)$. Hence, $[n \circ f]=[g]$. 
    
    We now check that $n$ induces the trivial variation. Observe that $n$ is the identity apart from on a local ball $B$ so given $\alpha \in H_2(M)$, there exists a representative $x$ of $\alpha$ disjoint from $B$. Therefore, 
    \[\Delta_n(\alpha)=[x-n(x)]=[x-x]=\Delta_{\id}(\alpha).\]
\end{proof}

We now prove the realization result in the smooth category. 

\begin{prop}\label[prop]{sm_surj}
    Let $M_0, M_1$ be homeomorphic simply-connected smooth $4$-manifolds. Fix a homeomorphism $F\maps M_1 \to M_0$ with $F\vert_{\del M_1}=h^{-1}$  for some $h$ and $\Delta \in \cV(H_2(M_0), \lambda_{M_0})$. Then, there exists a smooth $h$-cobordism $(W,M_0, f_0, M_1, f_1,g)$ rel. $h$ such that $\Delta_{F\circ \varphi_W}=\Delta$. Furthermore, if $M_0$ is spin, then we can also choose $W$ so that $\Theta(\varphi_W)=x$ for any $x\in H^1(M_0, \del M_0;\ZZ/2)$.
\end{prop}
\begin{proof}
    Let $\Delta\in \cV(H_2(M_0), \lambda_{M_0})$. As in the topological setting, \cite[Theorem 3.8]{OrsonPowell} provides a homeomorphism $G\maps M_0\to M_0$ restricting to the identity on the boundary such that $\Delta_G=\Delta$. By \cref{pinchy}, there exists a homotopy equivalence $n\maps M_0 \to M_0$ inducing the trivial variation such that \[[n \circ F]=[G]\in \cS^{\Diff}(M_0, \del M_0),\] 
    so by \cref{ses}, there exists a smooth $h$-cobordism $(W, M_0, f_0, M_1, f_1,g)$ from $n\circ F$ to $G$. 
    
    We need to check that $W$ induces the correct variation. Recall from diagram \ref{equiv} in the definition of the structure set that there exists a homotopy equivalence $\Phi$ restricting to a diffeomorphism on $\del_V W$ such that the following diagram commutes. 
    \[\begin{tikzcd}
M_1 \arrow[r, "n \circ F"] \arrow[d, "f_1"] & M_0 \arrow[d, "\id"]                    \\
\del_+(W) \arrow[r] \arrow[d, "i_1"', hook] & M_0\times \{1\} \arrow[d, "j_1", hook'] \\
W \arrow[r, "\Phi"]                                 & M_0 \times I                            \\
\del_-(W) \arrow[u, "i_-", hook'] \arrow[r] & M_0\times \{0\} \arrow[u, "j_0"', hook] \\
M_0 \arrow[r, "G"] \arrow[u, "f_0"]         & M_0. \arrow[u, "\id"']                  
\end{tikzcd}\]
    Hence, $\varphi_W\simeq F^{-1} \circ n^{-1}\circ G$ and thus, $\Delta_{F\circ \varphi_W}=\Delta_{n\circ G}=\Delta_{\id}*\Delta_G=\Delta$ by \cref{htpy} and \cref{var_hom}.
    
    Note that $W$ may not be an $h$-cobordism rel. $h$. To fix this, recall from \cref{replace_g} that we have some flexibility in how we choose $g$. Namely, we can replace $g$ with any other homeomorphism $g'$ such that $f_0^{-1}\circ g'_0=\id$. In particular, if we choose $g'$ to be $(\Phi\vert_{\del_V W})^{-1}$ (noting that $G\vert_{\del M_0}=\id$ so $f_0^{-1}\circ g'_0=\id$ as needed), then we can restrict the above diagram to obtain
    \[\begin{tikzcd}
\del M_1 \arrow[r, "n \circ F"] \arrow[d, "f_1"] & \del M_0 \arrow[d, "\id"]               \\
\del\del_+(W)                                    & \del M_0\times \{1\}, \arrow[l, "g_1'"']
\end{tikzcd}\]
where every map is a homeomorphism. So, recalling that $n\vert_{\del M_0}=\id$, we have $n\circ F\vert_{\del M_1}=F\vert_{\del M_1}$ so $f_1^{-1}\circ g_1'=(n\circ F)^{-1}=h^{-1}$ on $\del M_0$ and we obtain an $h$-cobordism rel. $h$. 

If $M_0$ is spin, we can choose $G$ from \cite[Theorem 3.8]{OrsonPowell} such that $\Theta(G)=x-\Theta(n^{-1})$ for any given $x\in H^1(M_0, \del M_0;\ZZ/2)$, so 
\[\Theta(F\circ \varphi_W)=\Theta(F\circ F^{-1}\circ n^{-1}\circ G)=\Theta(\id)+\Theta(n^{-1})+\Theta(G)=x.\]
\end{proof}

\noindent\textit{Proof of \cref{surj}}. This follows immediately from \cref{top_surj} and \cref{sm_surj}. \qed

\section{Uniqueness of $h$-cobordisms}\label[section]{uniqueness}
The goal of this section is to prove the injectivity portion of the main theorem:
\begin{thm}\label[thm]{inj}
Let $M_0, M_1$ be homeomorphic simply-connected $\CAT$ $4$-manifolds and fix a homeomorphism $F\maps M_1 \to M_0$ with $F\vert_{\del M_1}=h^{-1}$ for some $h$. Let $(W, M_0, f_0, M_1, f_1,g)$, $(W', M_0, f_0', M_1, f_1',g')$ be $h$-cobordisms rel. $h$ such that $\Delta_{F\circ \varphi_{W}}=\Delta_{F\circ \varphi_{W'}}$. If the $M_i$ are spin, assume additionally that $\Theta(F\circ \varphi_{W})=\Theta(F\circ \varphi_{W'})$. Then, $W$ is $\CAT$-isomorphic rel. boundary to $W'$. 
\end{thm}

Throughout this section, denote $X:=W \cup_{\del W} -W'$ glued via the map \[(f_0'\circ f_0^{-1})\cup(g'\circ g^{-1})\cup(f_1'\circ f_1^{-1})\maps \del W \to \del W'.\] We need to prove the following items:
\begin{enumerate}
    \item There exist $(B,\xi)$-structures on $M_0,M_1$ extending over $X$.
    \item $X$ is $(B,\xi)$-nullbordant.
    \item The surgery obstruction is elementary. 
\end{enumerate}
Recall that in the spin setting, $B=\BCATSpin$. Hence, we prove the following lemma. 
\begin{lemma}\label[lemma]{Xspin}
   Let $M_0, M_1$ be homeomorphic simply-connected spin $\CAT$ $4$-manifolds. Fix a homeomorphism $F\maps M_1 \to M_0$ with $F\vert_{\del M_1}=h^{-1}$ for some $h$ and let $(W, M_0, f_0, M_1, f_1,g)$, $(W', M_0, f_0', M_1, f_1',g')$ be $h$-cobordisms rel. $h$ such that $\Theta(F\circ \varphi_{W})=\Theta(F\circ \varphi_{W'})$. Then, $X$ is spin. 
\end{lemma}

\begin{proof}
    We will construct a spin structure on $X$ by finding spin structures on $\del W,\del W'$ that map to each other under the gluing map and extend to spin structures on $W,W'$ respectively. 

    Let $\fs$ be a spin structure on $M_0$ and denote $\fs_V$ to be the product spin structure coming from the restriction of $\fs$ to the boundary on $\del M_0 \times [0,1]$. Observe that \[(f_0^{-1})^{*}(\fs) \cup (g^{-1})^{*}(\fs_V)\cup(f_0^{-1}\circ i_-^{-1} \circ i_+^{-1})^*(\fs),\] 
    is a spin structure on $\del W$ that extends to one on $W$ (namely $((i_- \circ f_0)^*)^{-1} (\fs)$). Under the gluing map $(f_0'\circ f_0^{-1})\cup(g'\circ g^{-1})\cup(f_1'\circ f_1^{-1})$, this is sent to the spin structure 
    \begin{align*}
     &(f_0^{-1}\circ f_0\circ f_0'^{-1})^{*}(\fs) \cup (g^{-1}\circ g\circ g'^{-1})^{*}(\fs_V)\cup(f_0^{-1}\circ i_-^{-1} \circ i_+^{-1}\circ f_1\circ f_1'^{-1})^*(\fs)\\
     &= (f_0'^{-1})^{*} (\fs) \cup (g'^{-1})^* (\fs_V)\cup (f_0^{-1}\circ i_-^{-1}\circ i_+^{-1}\circ f_1\circ f_1'^{-1})^* (\fs),
    \end{align*}
    on $\del W'$. We need to show that this extends to a spin structure on $W'$. Note that \[(f_0'^{-1})^* (\fs) \cup (g'^{-1})^*(\fs_V)\cup (f_0'^{-1}\circ i_-'^{-1}\circ i_+')^* (\fs),\] 
    is a spin structure on $\del W'$ that extends. Because $\Theta(F\circ\varphi_W)=\Theta(F\circ\varphi_{W'})$, we have by \cite[Lemma 2.7]{OrsonPowell} that 
    \[\Theta(\varphi_W^{-1}\circ \varphi_{W'})=\Theta(\varphi_W^{-1}\circ F^{-1}\circ F\circ \varphi_{W'})=-\Theta(F\circ\varphi_W)+ \Theta( F\circ \varphi_{W'})=0=\Theta(\id),\]
    so $(\varphi_{W}^{-1}\circ \varphi_{W'})^*(\fs)=\fs$. Hence, 
    \begin{align*}
    (f_0'^{-1}\circ i_-'^{-1}\circ i_+')^*(\fs)&=(\varphi_W^{-1}\circ \varphi_W'\circ f_0'^{-1}\circ i_-'^{-1}\circ i_+')^* (\fs)\\
    &=(f_0^{-1}\circ i_-^{-1}\circ i_+ \circ f_1 \circ f_1'^{-1}\circ i_+'^{-1}\circ i_-' \circ f_0' \circ f_0'^{-1}\circ i_-'^{-1}\circ i_+')^* (\fs)\\
    &=(f_0^{-1}\circ i_-^{-1}\circ i_+^{-1}\circ f_1\circ f_1'^{-1})^* (\fs),
    \end{align*}
    which implies that the spin structure from the gluing extends over $W'$.
\end{proof}
For the next two steps, it will be convenient to know the homology groups of $X$ (in particular, due to \cref{rk_elem}).

\begin{lemma}\label[lemma]{H grps}
Let $M_0, M_1$ be homeomorphic simply-connected $\CAT$ $4$-manifolds. Fix a homeomorphism $F\maps M_1 \to M_0$ with $F\vert_{\del M_1}=h^{-1}$ for some $h$ and let $(W, M_0, f_0, M_1, f_1,g)$, $(W', M_0, f_0', M_1, f_1',g')$ be $h$-cobordisms rel. $h$ such that $\Delta_{F\circ \varphi_{W}}=\Delta_{F\circ \varphi_{W'}}$. Then, there exists an isomorphism \[H_3(X) \cong H_2(M_0,\del M_0)\oplus H_3(W).\] 
\end{lemma}
\begin{proof}
To simplify notation, let $i_0=i_- \circ f_0$, $i_1=i_+\circ f_1$, and likewise for $W'$. We also write $\bar{f}_*$ for induced maps on relative homology groups to distinguish from those on absolute homology. Lastly, we note that there exist homeomorphisms such that the following square commutes and hence, we write $M_0 \cup_h -M_1$ everywhere instead of $M_0 \cup_{\id} \del M_0 \times [0,1] \cup_h -M_1$:
\[\begin{tikzcd}
{M_0 \cup_{\id} \del M_0\times [0,1] \cup_h -M_1} \arrow[r, "f_0\cup g \cup f_1"] \arrow[d, "\cong"] & \del_- W \cup_{\id} \del_V W \cup \del_+ W \arrow[d, "\cong"] \\
M_0 \cup_h -M_1 \arrow[r, "f_0 \cup f_1"]                                                       & \del_- W \cup \del_+ W                                  
.\end{tikzcd}\]

Let $j,j'$ be the inclusions of $M_0 \cup M_1$ into $W,W'$ respectively. The Mayer Vietoris sequence for $X=W \cup W'$ gives us

\[ H_3(M_0 \cup M_1) \xrightarrow{(j_*, j_*')} H_3(W)\oplus H_3(W') \ra H_3(X)\to H_2(M_0 \cup M_1)\xrightarrow{(j_*,j_*')}H_2(W)\oplus H_2(W).\]

Observe that the isomorphism ${i_{0}}_{*} \maps H_3(M_0)\ra H_3(W)$ factors through $H_3(M_0 \cup_h -M_1)$ and the map $j_*\maps H_3(M_0)\ra H_3(M_0 \cup_h -M_1)$ is surjective because  \[H_3(M_0 \cup_h -M_1, M_0)\cong H_3(M_1, \del M_1) \cong 0,\] by excision. Hence, the inclusion induces a surjection $j_* \maps H_3(M_0 \cup_h -M_1) \twoheadrightarrow H_3(W)$. Noting that \[H_4(W, M_0 \cup_h -M_1)\cong H^1(W)\cong 0,\] 
we obtain an isomorphism $H_3(M_0 \cup_h -M_1)\cong H_3(W)$ and likewise for $W'$. In particular, this implies that the image of $(j_*, j_*')\maps H_3(M_0 \cup M_1) \to H_3(W)\oplus H_3(W')$ is isomorphic $H_3(W)$ because if we identify $H_3(W')$ with $H_3(W)$ via $j_*\circ j_*'^{-1}$, then the image is the diagonal of $H_3(W)\oplus H_3(W)$. Therefore, we obtain a short exact sequence
\[0 \ra H_3(W) \ra H_3(X) \ra \ker(j_*, j'_*) \ra 0.\]

For the remainder of this proof, we use $j_*,j_*'$ only to refer to the maps on $H_2(M_0\cup H_1)$. Observe that it will suffice to show that $\ker j_* \cong H_2(M_0, \del M_0)$ and $\ker j_* = \ker j'_*$. To see this, note that if $\ker j_* = \ker j'_*$, then $\ker(j_*, j'_*)=\ker j_*$, and if $\ker j_*\cong H_2(M_0, \del M_0)$, then the short exact sequence splits (because $H_2(M_0, \del M_0)$ is free) to give us $H_3(X) \cong H_2(M_0,\del M_0)\oplus H_3(W)$ as desired. 

By Mayer-Vietoris, \[H_3(M_0 \cup_h -M_1,\del M_0) \cong H_3(M_0, \del M_0)\oplus H_3(M_1, \del M_1) \cong 0,\] 
and recalling that $W\simeq M_0$ rel. $\del M_0$, we get \[H_3(W,\del M_0) \cong H_3(M_0, \del M_0)\cong 0.\] 
Hence, the inclusion induces the following commutative diagram between the long exact sequences of pairs:
\[\begin{tikzcd}
0 \arrow[r] \arrow[d] & H_2(\del M_0) \arrow[r] \arrow[d, "\id"] & H_2(M_0 \cup_h -M_1) \arrow[r, "\psi"] \arrow[d, "j_*"] & {H_2(M_0 \cup_h -M_1, \del M_0)} \arrow[d, "\bar j_*"] \\
0 \arrow[r]           & H_2(\del M_0) \arrow[r]           & H_2(W) \arrow[r, "\psi_W"]                         & {H_2(W, \del M_0)}.                                 
\end{tikzcd}\]
We first show that $\psi$ gives an isomorphism $\ker j_* \cong \psi(\ker j_*)$. Surjectivity is trivial so let $x\in \ker j_* \cap \ker \psi$. By exactness, there exists $y\in \ker H_2(\del M_0)$ mapping to $x$. But by commutativity, $y \mapsto j_*(x)=0 \in H_2(W)$ so $y=0$ and hence, $x=0$. 

Next, we show that $\psi(\ker j_*)=\ker \bar {j_*}$. 
The long exact sequence of a triple fits into the following commutative diagram with the long exact sequence of a pair, with $\psi, \psi_W$ induced by the inclusion:
\[\begin{tikzcd}
{H_3(W, M_0 \cup_h -M_1)} \arrow[r, "\del_W"] \arrow[d, "\id"] & H_2(M_0 \cup_h -M_1) \arrow[r, "j_*"] \arrow[d, "\psi"] & H_2(W) \arrow[d, "\psi_{W}"] \\
{H_3(W, M_0 \cup_h -M_1)} \arrow[r, "\del"]            & {M_2(M_0 \cup_h -M_1, \del M_0)} \arrow[r, "\bar j_*"]     & {H_2(W, \del M_0)}.              
\end{tikzcd}\]
Given $\psi(x)\in \psi(\ker j_*)$, we have by commutativity $\bar j_* (\psi(x))=\psi_W \circ j_* (x)=0$ so $\psi(\ker j_*)\subset \ker \bar {j_*}$. For the reverse inclusion, let $x \in \ker \bar j_*$. By exactness, there exists $y \in $ $H_3(W, M_0 \cup_h -M_1)$ such that $\del(y)=x$. By exactness, $\del_W(y)\in \ker j_*$ and by commutativity $\psi(\del_W y)=x$. Hence, $x \in \psi(\ker j_*))$ as desired.

Observe that that by Mayer-Vietoris, we have an isomorphism \[H_2(M_0 \cup_h -M_1, \del M_0) \to H_2(M_0, \del M_0)\oplus H_2(M_1, \del M_1),\] given by $([x_0],[x_1])\mapsto [x_0+x_1]$. Composing with $\bar j_*$, this gives the map 
\[\bar{i_0}_*+\bar{i_1}_*\maps H_2(M_0,\del M_0)\oplus H_2(M_1,\del M_1)\to H_2(W, \del M_0).\]
Then, recalling that $i_0 \maps M_0 \to W$ is a homotopy equivalence, 
\[\ker (\bar{i_0}_*+ \bar{i_1}_*)=\ker(\bar{i_0}_*^{-1}\circ(\bar{i_0}_*+ \bar{i_1}_*))=\ker(\id+ \bar{i_0}_*^{-1}\circ\bar{i_1}_*))=\Graph(-\bar{i_0}_*^{-1}\circ\bar{i_1}_*).\]
Putting everything together, we get 
\[\ker j_*\cong \psi(\ker j_*)=\ker \bar j_* \cong \ker (\bar{i_0}_*+\bar{i_1}_*)\cong \Graph(-\bar{(i_0)}_*^{-1}\circ\bar{(i_1)}_*)\cong H_2(M_0, \del M_0).\] 

Now, we show that $\ker j_* \subset \ker j'_*$.  Let $[x]\in\ker j_*$. Again, via the isomorphism \[H_2(M_0, \del M_0)\oplus H_2(M_1, \del M_1)\cong H_2(M_0 \cup_h -M_1, \del M_0),\] 
we can write $\psi([x])=[x_0 +x_1$], where $x_0, x_1$ are relative chains in $C_2(M_0, \del M_0)$ and $C_2(M_1, \del M_1)$ respectively. Hence, $[x]=[x_0+x_1+x_\del]$, where by an abuse of notation, we regard $x_0 \in C_2(M_0)$, $x_1 \in C_2(M_1)$, and $x_\del \in C_2(\del M_0)$. We can do this because $\psi$ is induced from the quotient map \[C_2(M_0\cup M_1)\to C_2(M_0\cup M_1)/C_2(M_0),\] 
so we can take $x_0$ and $x_1$ to be representatives for their respective cosets in the quotient. Observe that $\Delta_{F\circ \varphi_W}=\Delta_{F\circ \varphi_{W'}}$ implies $\Delta_{\varphi_{W'}^{-1}\circ \varphi_{W}}=0$ by \cref{var_hom}. Hence, we have 
\begin{align*}
    (i_0'^{-1}\circ j')_*([x])
    &=(i_0'^{-1}\circ j')_*([x])-(i_0^{-1}\circ j)_*([x])\\
    &=[i_0'^{-1}\circ i_0'(x_0)+\varphi_{W'}^{-1}(x_1)+i_0'^{-1}\circ i_0'(x_\del)-i_0^{-1}\circ i_0(x_0)-\varphi_{W}^{-1}(x_1)+i_0^{-1}\circ i_0(x_\del)]\\
    &=[\varphi_{W'}^{-1}(x_1)-\varphi_{W}^{-1}(x_1)]\\
    &=[\varphi_{W'}^{-1}(x_1)-\varphi_{W}^{-1}(x_1)]-\Delta_{\varphi_{W}^{-1}\circ \varphi_{W'}}([\varphi_{W'}^{-1}(x_1)])\\
    &=[\varphi_{W'}^{-1}(x_1)-\varphi_{W}^{-1}(x_1)-\varphi_{W'}^{-1}(x_1)+\varphi_{W}^{-1}\circ \varphi_{W'} \circ \varphi_{W'}^{-1}(x_1)]\\
    &=0
\end{align*}

In the second equality, we used that by $j'$ was the inclusion so it restricts to $i_0'$ on $M_0$ and $i_1'$ on $M_1$. Noting that $i_0'$ is an isomorphism, we obtain that $j'_*([x])=0$ and $\ker j_* \subset \ker j_*'$. An identical argument provides the other inclusion so we have equality. 
\end{proof}

\begin{cor}\label[lemma]{H_cor}
    The homology groups of $X$ are 
    \[H_n(X)=\begin{cases}
        \ZZ & n=0,5\\
        0&n=1,4 \\
        \ZZ^{b_0(\del M_0)-1}\oplus \ZZ^{b_2(M_0)} &n=2,3.
    \end{cases}\]
\end{cor}
\begin{proof}
$X$ is simply connected because $W,W'$ were simply connected and $\del W$ is connected. Also, $X$ is a closed $5$-manifold so $H_5(X)=\ZZ$. Lastly, Poincare duality, universal coefficients, and \cref{H grps} give \[H_2(X)\cong H_3(X)^*\cong H_2(M_0,\del M_0)^*\oplus H_3(W)^*\cong\ZZ^{b_0(\del M_0)-1}\oplus \ZZ^{b_2(M_0)}\]
and 
$H_4(X)\cong H_1(X)^*=0$. 
\end{proof}
We will now apply \cref{rk_elem} to show that $W\cong W'$. 

\noindent\textit{Proof of \cref{inj}.} Recall that the normal $1$-type of a simply connected $4$-manifold is equivalent to a choice of orientation in the nonspin case and a choice of spin structure in the spin case. By definition of an $h$-cobordism, an orientation on $\del W$ must extend over $W$. In the spin case, we have by \cref{Xspin} that we can choose spin structures on $\del W$ extending to a spin structure on $X$. Hence, there exist normal $1$-smoothings on $W$, $W'$ that extend to $X$. 

Now, we seek to show that $X$ is $(B,\xi)$-nullbordant. In the spin case, $\Omega_5^{\CATSpin}\cong0$. In the nonspin case, $\Omega_5^{\SCAT}\cong\ZZ/2$ and bordism classes are detected by $w_2 \cdot w_3$. By \cite{LMP}, $w_2\cdot w_3$ is equal to $\sigma(X;\ZZ/2)-\sigma(X;\QQ)$ where $\sigma(X;F):=\sum_i\rk(H_i(X;F))$ for a field $F$. By \cref{H_cor}, $H_i(X)$ is free for all $i$ so by the UCT, $\sigma(X;\QQ)-\sigma(X;\ZZ/2)=0$. Hence, in both scenarios, $X$ is $(B,\xi)$-nullbordant via some nullbordism $Y$, which we can surger so that $Y\to B$ is $3$-connected. 

We want to use \cref{rk_elem} to show that the surgery obstruction is elementary. From the calculation of $H_*(X)$ in \cref{H_cor}, we see that $H_2(W), H_2(W')$ are free. Also, 
\[\ker(H_3(Y)\ra H_3(Y,W))=\im(H_3(W)\ra H_3(Y))=\im(H_3(M_0)\ra H_3(Y)),\] 
and likewise for $W'$ because the inclusion $M_0\hookrightarrow Y$ factors through $W,W'$ via homotopy equivalences so \cref{rk_elem} is applicable. From \cref{H_cor}, we also see that $\rk H_3(X)=\rk H_2(W)+ \rk H_3(W)$. 
Hence, \cref{rk_elem} shows that the surgery obstruction is elementary and as a result, $W \cong W'$ rel. boundary, which is the conclusion of \cref{inj}. 
\qed

\section{Proof of the Main Theorem}\label[section]{proof}
We now compile our results and prove the main theorem. 

\noindent\textit{Proof of \cref{thm}.} 
By \cref{inj} and \cref{surj}, the maps \[\cH^{\CAT}(M_0,M_1,h)\to \cV(H_2(M_0),\lambda_{M_0}),\] 
and \[\cH^{\CAT}(M_0,M_1,h)\to H^1(M_0,\del M_0)\times \cV(H_2(M_0),\lambda_{M_0}),\] 
given by $W \mapsto \Delta_{F\circ \varphi_{W}}$ and $W \mapsto (\Delta_{F\circ \varphi_{W}},\Theta(F\circ \varphi_{W}))$ are bijective. In the case that $M_0=M_1=\id$, we choose $F=\id$ and these maps are group homomorphisms via \cref{hcob_var_hom} and \cref{hcob_spin_hom}.
\qed
\section{Pseudoisotopies}\label[section]{pseudoisotopies}
We conclude by recalling the definition of a pseudoisotopy and proving injectivity of the pseudo-mapping class group computation in \cite{OrsonPowell}. 
\begin{defn}
    Let $F\maps M\to M$ be a homeomorphism with $F\vert_{\del M}=\id$. We say that $F$ is \emph{pseudoisotopic} to the identity if there exists a homeomorphism $\Phi \maps M\times [0,1]  \to M \times [0,1]$ such that $\Phi\vert_{M\times\{0\}\cup \del M \times [0,1]}=\id$ and $\Phi\vert_{\del M\times\{1\}}=F$. We call $\Phi$ a \emph{pseudoisotopy} and denote as $\tilde{\pi}_0\Homeo^+(M,\partial M)$ the group of orientation preserving homeomorphisms of $M$ restricting to the identity on $\del M$ modulo pseudoisotopy. The same definitions also all hold in the smooth category.
\end{defn}

Observe that a pseudoisotopy $\Phi$ from $F$ to the identity is exactly the same as a $\CAT$-isomorphism rel. boundary between the $h$-cobordisms $(M\times[0,1], M, \id, M, F^{-1}, \id)$ and $(M\times[0,1], M, \id, M, \id, \id)$. Hence, we obtain the following alternate proof of \cite[Theorem 3.1]{OrsonPowell}.
\begin{thm}[{\cite[Theorem 3.1]{OrsonPowell}}]
    Let $M$ be a compact, simply connected $\CAT$ $4$-manifold.
\begin{enumerate}
\item When $X$ is spin, the homomorphisms
\[
\begin{array}{ccl}
\tilde{\pi}_0\Homeo^+(M,\partial M)&\hookrightarrow&H^1(M,\partial M;\Z/2)\times\cV(H_2(M),\lambda_M), \quad{and}\\
\tilde{\pi}_0\Diffeo^+(M,\partial M)&\hookrightarrow&H^1(M,\partial M;\Z/2)\times\cV(H_2(M),\lambda_M),
\end{array}
\]
given by $F\mapsto (\Theta(F),\Delta_F)$, are injective.
\item When $X$ is not spin, the homomorphisms
\[
\begin{array}{ccl}
\widetilde{\pi}_0\Homeo^+(M,\partial M)&\hookrightarrow&\mathcal{V}(H_2(M),\lambda_M), \quad{and}\\
\widetilde{\pi}_0\Diffeo^+(M,\partial M)&\hookrightarrow&\mathcal{V}(H_2(M),\lambda_M),
\end{array}
\]
given by $F\mapsto \Delta_F$, are injective.
\end{enumerate}
\end{thm}
\begin{proof}
Let $F\maps M \to M$ be a $\CAT$-isomorphism such that $\Delta_F$ is trivial (and $\Theta(F)$ if $M$ is spin). Denoting $W,W'$ to be the $h$-cobordisms $(M\times[0,1], M, \id, M, F^{-1}, \id)$ and $(M\times[0,1], M, \id, M, \id, \id)$ respectively, we have $\varphi_W\simeq F$ and $\varphi_{W'}\simeq \id$. In particular, $\Delta_{\varphi_{W}}=\Delta_{F}=0$, and likewise for $\Theta$ if $M$ is spin. Hence, by \cref{thm}, $W\cong W'$ rel. boundary and hence, by the above observation, there exists a pseudoisotopy taking $F$ to the identity. 
\end{proof}

As a final remark, we note that, as in \cite{OrsonPowell}, we can upgrade this pseudoisotopy computation to the level of isotopy via work of Quinn and Perron (\cite{Quinn,Perron}).

\bibliographystyle{alpha}
\bibliography{references}

@article {Saeki,
    AUTHOR = {Saeki, Osamu},
     TITLE = {Stable mapping class groups of 4-manifolds with boundary},
   JOURNAL = {Trans. Amer. Math. Soc.},
  FJOURNAL = {Transactions of the American Mathematical Society},
    VOLUME = {358},
      YEAR = {2006},
    NUMBER = {5},
     PAGES = {2091--2104},
      ISSN = {0002-9947},
   MRCLASS = {57R52 (57N13 57N37)},
  MRNUMBER = {2197449},
MRREVIEWER = {Raphael Zentner},
       DOI = {10.1090/S0002-9947-05-03748-7},
       URL = {https://doi.org/10.1090/S0002-9947-05-03748-7},
}

@article {OrsonPowell,
    AUTHOR = {Orson, Patrick and Powell, Mark},
     TITLE = {Mapping class groups of simply connected 4-manifolds with
              boundary},
   JOURNAL = {J. Differential Geom.},
  FJOURNAL = {Journal of Differential Geometry},
    VOLUME = {131},
      YEAR = {2025},
    NUMBER = {1},
     PAGES = {199--275},
      ISSN = {0022-040X},
   MRCLASS = {57K40 (57N37 57R52 57S05)},
  MRNUMBER = {4947553},
       DOI = {10.4310/jdg/1755544135},
       URL = {https://doi.org/10.4310/jdg/1755544135},
}

@book {Wall,
    AUTHOR = {Wall, C. T. C.},
     TITLE = {Surgery on compact manifolds},
    SERIES = {Mathematical Surveys and Monographs},
    VOLUME = {69},
   EDITION = {Second},
      NOTE = {Edited and with a foreword by A. A. Ranicki},
 PUBLISHER = {American Mathematical Society, Providence, RI},
      YEAR = {1999},
     PAGES = {xvi+302},
      ISBN = {0-8218-0942-3},
   MRCLASS = {57R67 (19J25 57-02)},
  MRNUMBER = {1687388},
       DOI = {10.1090/surv/069},
       URL = {https://doi.org/10.1090/surv/069},
}

@article {CochranHabegger,
    AUTHOR = {Cochran, Tim D. and Habegger, Nathan},
     TITLE = {On the homotopy theory of simply connected four manifolds},
   JOURNAL = {Topology},
  FJOURNAL = {Topology. An International Journal of Mathematics},
    VOLUME = {29},
      YEAR = {1990},
    NUMBER = {4},
     PAGES = {419--440},
      ISSN = {0040-9383},
   MRCLASS = {57N13 (55P10 55Q52)},
  MRNUMBER = {1071367},
MRREVIEWER = {Laurence R. Taylor},
       DOI = {10.1016/0040-9383(90)90014-B},
       URL = {https://doi.org/10.1016/0040-9383(90)90014-B},
}

@article {LMP,
    AUTHOR = {Lusztig, G. and Milnor, J. and Peterson, F. P.},
     TITLE = {Semi-characteristics and cobordism},
   JOURNAL = {Topology},
  FJOURNAL = {Topology. An International Journal of Mathematics},
    VOLUME = {8},
      YEAR = {1969},
     PAGES = {357--359},
      ISSN = {0040-9383},
   MRCLASS = {57.10},
  MRNUMBER = {246308},
MRREVIEWER = {H. Suzuki},
       DOI = {10.1016/0040-9383(69)90021-4},
       URL = {https://doi.org/10.1016/0040-9383(69)90021-4},
}

@article {SAD,
    AUTHOR = {Kreck, Matthias},
     TITLE = {Surgery and duality},
   JOURNAL = {Ann. of Math. (2)},
  FJOURNAL = {Annals of Mathematics. Second Series},
    VOLUME = {149},
      YEAR = {1999},
    NUMBER = {3},
     PAGES = {707--754},
      ISSN = {0003-486X},
   MRCLASS = {57R67 (57R65)},
  MRNUMBER = {1709301},
MRREVIEWER = {Laurence R. Taylor},
       DOI = {10.2307/121071},
       URL = {https://doi.org/10.2307/121071},
}

@article {SixtMonoids,
    AUTHOR = {Sixt, J\"{o}rg},
     TITLE = {Even-dimensional {$l$}-monoids and {$L$}-theory},
   JOURNAL = {Algebr. Geom. Topol.},
  FJOURNAL = {Algebraic \& Geometric Topology},
    VOLUME = {7},
      YEAR = {2007},
     PAGES = {479--515},
      ISSN = {1472-2747},
   MRCLASS = {57R67 (19G24 57R65)},
  MRNUMBER = {2308954},
MRREVIEWER = {A. A. Ranicki},
       DOI = {10.2140/agt.2007.7.479},
       URL = {https://doi.org/10.2140/agt.2007.7.479},
}

@article {KPLT,
    AUTHOR = {Kasprowski, Daniel and Land, Markus and Powell, Mark and
              Teichner, Peter},
     TITLE = {Stable classification of 4-manifolds with 3-manifold
              fundamental groups},
   JOURNAL = {J. Topol.},
  FJOURNAL = {Journal of Topology},
    VOLUME = {10},
      YEAR = {2017},
    NUMBER = {3},
     PAGES = {827--881},
      ISSN = {1753-8416},
   MRCLASS = {57N13 (57N70)},
  MRNUMBER = {3797598},
MRREVIEWER = {Richard Stong},
       DOI = {10.1112/topo.12025},
       URL = {https://doi.org/10.1112/topo.12025},
}

@article {Kreck_hcob,
    AUTHOR = {Kreck, Matthias},
     TITLE = {{$h$}-cobordisms between 1-connected 4-manifolds},
   JOURNAL = {Geom. Topol.},
  FJOURNAL = {Geometry and Topology},
    VOLUME = {5},
      YEAR = {2001},
     PAGES = {1--6},
      ISSN = {1465-3060,1364-0380},
   MRCLASS = {57R80 (57N13 57Q20)},
  MRNUMBER = {1812433},
MRREVIEWER = {Terry\ Lawson},
       DOI = {10.2140/gt.2001.5.1},
       URL = {https://doi.org/10.2140/gt.2001.5.1},
}

@book{Hatcher,
	address = {New York},
	title = {Algebraic topology},
	isbn = {978-0-521-79540-1},
	language = {eng},
	publisher = {Cambridge University Press},
	author = {Hatcher, Allen},
	year = {2001},
}

@misc{Hatchervb,
  author = {Hatcher, Allen},
  title = {Vector Bundles and K-Theory},
  year = {2003},
  howpublished={\url{https://pi.math.cornell.edu/~hatcher/VBKT/VB.pdf}}
       }

@article {Boyer,
    AUTHOR = {Boyer, Steven},
     TITLE = {Simply-connected {$4$}-manifolds with a given boundary},
   JOURNAL = {Trans. Amer. Math. Soc.},
  FJOURNAL = {Transactions of the American Mathematical Society},
    VOLUME = {298},
      YEAR = {1986},
    NUMBER = {1},
     PAGES = {331--357},
      ISSN = {0002-9947,1088-6850},
   MRCLASS = {57N13 (57N10)},
  MRNUMBER = {857447},
MRREVIEWER = {V.\ G.\ Turaev},
       DOI = {10.2307/2000623},
       URL = {https://doi.org/10.2307/2000623},
}

@book {Brow72,
    AUTHOR = {Browder, William},
     TITLE = {Surgery on simply-connected manifolds},
    SERIES = {Ergebnisse der Mathematik und ihrer Grenzgebiete [Results in
              Mathematics and Related Areas], Band 65},
 PUBLISHER = {Springer-Verlag, New York-Heidelberg},
      YEAR = {1972},
     PAGES = {ix+132},
   MRCLASS = {57D65},
  MRNUMBER = {358813},
MRREVIEWER = {P.\ J.\ Kahn},
}

@incollection {Brow93,
    AUTHOR = {Browder, William},
     TITLE = {Homotopy type of differentiable manifolds},
 BOOKTITLE = {Novikov conjectures, index theorems and rigidity, {V}ol. 1
              ({O}berwolfach, 1993)},
    SERIES = {London Math. Soc. Lecture Note Ser.},
    VOLUME = {226},
     PAGES = {97--100},
 PUBLISHER = {Cambridge Univ. Press, Cambridge},
      YEAR = {1995},
      ISBN = {0-521-49796-5},
   MRCLASS = {57R65 (57R55)},
  MRNUMBER = {1388298},
MRREVIEWER = {Terry\ Lawson},
       DOI = {10.1017/CBO9780511662676.006},
       URL = {https://doi.org/10.1017/CBO9780511662676.006},
}

@article {Nov,
    AUTHOR = {Novikov, S. P.},
     TITLE = {Homotopically equivalent smooth manifolds. {I}},
   JOURNAL = {Izv. Akad. Nauk SSSR Ser. Mat.},
  FJOURNAL = {Izvestiya Akademii Nauk SSSR. Seriya Matematicheskaya},
    VOLUME = {28},
      YEAR = {1964},
     PAGES = {365--474},
      ISSN = {0373-2436},
   MRCLASS = {57.10},
  MRNUMBER = {162246},
MRREVIEWER = {A.\ H.\ Wallace},
}

@book {Sull,
    AUTHOR = {Sullivan, Dennis Parnell},
     TITLE = {Triangulating homotopy equivalences},
      NOTE = {Thesis (Ph.D.)--Princeton University},
 PUBLISHER = {ProQuest LLC, Ann Arbor, MI},
      YEAR = {1966},
     PAGES = {200},
   MRCLASS = {Thesis},
  MRNUMBER = {2615489},
       URL =
              {http://gateway.proquest.com/openurl?url_ver=Z39.88-2004&rft_val_fmt=info:ofi/fmt:kev:mtx:dissertation&res_dat=xri:pqdiss&rft_dat=xri:pqdiss:6609643},
}

@book {Ranicki,
    AUTHOR = {Ranicki, Andrew},
     TITLE = {Algebraic and geometric surgery},
    SERIES = {Oxford Mathematical Monographs},
      NOTE = {Oxford Science Publications},
 PUBLISHER = {The Clarendon Press, Oxford University Press, Oxford},
      YEAR = {2002},
     PAGES = {xii+373},
      ISBN = {0-19-850924-3},
   MRCLASS = {57R65 (19J25 57R67)},
  MRNUMBER = {2061749},
MRREVIEWER = {Daniel\ Ruberman},
       DOI = {10.1093/acprof:oso/9780198509240.001.0001},
       URL = {https://doi.org/10.1093/acprof:oso/9780198509240.001.0001},
}

@incollection {Luck,
    AUTHOR = {L\"{u}ck, Wolfgang},
     TITLE = {A basic introduction to surgery theory},
 BOOKTITLE = {Topology of high-dimensional manifolds, {N}o. 1, 2 ({T}rieste,
              2001)},
    SERIES = {ICTP Lect. Notes},
    VOLUME = {9},
     PAGES = {1--224},
 PUBLISHER = {Abdus Salam Int. Cent. Theoret. Phys., Trieste},
      YEAR = {2002},
      ISBN = {92-95003-12-8},
   MRCLASS = {57R67 (57R65)},
  MRNUMBER = {1937016},
}

@book {FreedmanQuinn,
    AUTHOR = {Freedman, Michael H. and Quinn, Frank},
     TITLE = {Topology of 4-manifolds},
    SERIES = {Princeton Mathematical Series},
    VOLUME = {39},
 PUBLISHER = {Princeton University Press, Princeton, NJ},
      YEAR = {1990},
     PAGES = {viii+259},
      ISBN = {0-691-08577-3},
   MRCLASS = {57N13 (57-02)},
  MRNUMBER = {1201584},
MRREVIEWER = {Ian\ Hambleton},
}

@article {Perron,
    AUTHOR = {Perron, B.},
     TITLE = {Pseudo-isotopies et isotopies en dimension quatre dans la
              cat\'{e}gorie topologique},
   JOURNAL = {Topology},
  FJOURNAL = {Topology. An International Journal of Mathematics},
    VOLUME = {25},
      YEAR = {1986},
    NUMBER = {4},
     PAGES = {381--397},
      ISSN = {0040-9383},
   MRCLASS = {57N37 (57N13)},
  MRNUMBER = {862426},
MRREVIEWER = {L.\ Guillou},
       DOI = {10.1016/0040-9383(86)90018-2},
       URL = {https://doi.org/10.1016/0040-9383(86)90018-2},
}

@article {Quinn,
    AUTHOR = {Quinn, Frank},
     TITLE = {Isotopy of {$4$}-manifolds},
   JOURNAL = {J. Differential Geom.},
  FJOURNAL = {Journal of Differential Geometry},
    VOLUME = {24},
      YEAR = {1986},
    NUMBER = {3},
     PAGES = {343--372},
      ISSN = {0022-040X,1945-743X},
   MRCLASS = {57N13 (57N37 57S05)},
  MRNUMBER = {868975},
MRREVIEWER = {Laurence\ R.\ Taylor},
       URL = {http://projecteuclid.org/euclid.jdg/1214440552},
}

@article {structuresets,
    AUTHOR = {Cencelj, Matija and Muranov, Yuri V. and Repov\v{s},
              Du\v{s}an},
     TITLE = {On structure sets of manifold pairs},
   JOURNAL = {Homology Homotopy Appl.},
  FJOURNAL = {Homology, Homotopy and Applications},
    VOLUME = {11},
      YEAR = {2009},
    NUMBER = {2},
     PAGES = {195--222},
      ISSN = {1532-0073,1532-0081},
   MRCLASS = {57R67 (19J25)},
  MRNUMBER = {2591918},
       DOI = {10.4310/hha.2009.v11.n2.a10},
       URL = {https://doi.org/10.4310/hha.2009.v11.n2.a10},
}

@book {ranickiexactseq,
    AUTHOR = {Ranicki, Andrew},
     TITLE = {Exact sequences in the algebraic theory of surgery},
    SERIES = {Mathematical Notes},
    VOLUME = {26},
 PUBLISHER = {Princeton University Press, Princeton, NJ; University of Tokyo
              Press, Tokyo},
      YEAR = {1981},
     PAGES = {xvii+864},
      ISBN = {0-691-08276-6},
   MRCLASS = {57R67 (18F25 57-02)},
  MRNUMBER = {620795},
MRREVIEWER = {F.\ T.\ Farrell},
}

@article {lawson,
    AUTHOR = {Lawson, Terry},
     TITLE = {{$h$}-cobordisms between simply connected {$4$}-manifolds},
   JOURNAL = {Topology Appl.},
  FJOURNAL = {Topology and its Applications},
    VOLUME = {28},
      YEAR = {1988},
    NUMBER = {1},
     PAGES = {75--82},
      ISSN = {0166-8641,1879-3207},
   MRCLASS = {57N13},
  MRNUMBER = {927283},
MRREVIEWER = {Ian\ Hambleton},
       DOI = {10.1016/0166-8641(88)90037-5},
       URL = {https://doi.org/10.1016/0166-8641(88)90037-5},
}

@article {hollingsworth_morgan,
    AUTHOR = {Hollingsworth, John G. and Morgan, John W.},
     TITLE = {Homotopy triangulations and topological manifolds},
      YEAR = {1970},
}

@inproceedings {Brumfiel,
    AUTHOR = {Brumfiel, Gregory W.},
     TITLE = {Homotopy equivalences of almost smooth manifolds},
 BOOKTITLE = {Algebraic topology ({P}roc. {S}ympos. {P}ure {M}ath., {V}ol.
              {XXII}, {U}niv. {W}isconsin, {M}adison, {W}is., 1970)},
    SERIES = {Proc. Sympos. Pure Math., Vol. XXII},
     PAGES = {73--79},
 PUBLISHER = {Amer. Math. Soc., Providence, RI},
      YEAR = {1971},
   MRCLASS = {57D90},
  MRNUMBER = {319217},
MRREVIEWER = {R.\ K.\ Lashof},
}

@article {hambletonkreck,
    AUTHOR = {Hambleton, Ian and Kreck, Matthias},
     TITLE = {Homotopy self-equivalences of 4-manifolds},
   JOURNAL = {Math. Z.},
  FJOURNAL = {Mathematische Zeitschrift},
    VOLUME = {248},
      YEAR = {2004},
    NUMBER = {1},
     PAGES = {147--172},
      ISSN = {0025-5874,1432-1823},
   MRCLASS = {57N13 (55P10 57R80)},
  MRNUMBER = {2092726},
MRREVIEWER = {Terry\ Lawson},
       DOI = {10.1007/s00209-004-0657-9},
       URL = {https://doi.org/10.1007/s00209-004-0657-9},
}

@inproceedings {Brown,
    AUTHOR = {Brown, Morton},
     TITLE = {Locally flat embeddings of topological manifolds},
 BOOKTITLE = {Topology of 3-manifolds and related topics ({P}roc. {T}he
              {U}niv. of {G}eorgia {I}nstitute, 1961)},
     PAGES = {83--91},
 PUBLISHER = {Prentice-Hall, Inc., Englewood Cliffs, NJ},
      YEAR = {1961},
   MRCLASS = {54.78},
  MRNUMBER = {158373},
MRREVIEWER = {O.\ G.\ Harrold},
}

@misc{CCPS,
      title={Stably diffeomorphic manifolds and the realisation of modified surgery obstructions}, 
      author={Anthony Conway and Diarmuid Crowley and Mark Powell and Joerg Sixt},
      year={2024},
      eprint={2109.05632},
      archivePrefix={arXiv},
      primaryClass={math.GT},
      howpublished={\url{https://arxiv.org/abs/2109.05632}}
}

@book {may,
    AUTHOR = {May, J. P.},
     TITLE = {A concise course in algebraic topology},
    SERIES = {Chicago Lectures in Mathematics},
 PUBLISHER = {University of Chicago Press, Chicago, IL},
      YEAR = {1999},
     PAGES = {x+243},
      ISBN = {0-226-51182-0; 0-226-51183-9},
   MRCLASS = {55-02 (18-02 57-02)},
  MRNUMBER = {1702278},
MRREVIEWER = {R.\ M.\ Vogt},
}
\end{document}